\documentclass[journal]{IEEEtran}
\ifCLASSINFOpdf
\else
\fi
\IEEEoverridecommandlockouts
\usepackage{epsfig}
\usepackage{color}
\usepackage{amssymb,latexsym,amsfonts,amsmath}
\usepackage{graphicx}
\usepackage{cite}
\usepackage{flushend}
\allowdisplaybreaks[4]

\newtheorem{theorem}{Theorem}
\newtheorem{lemma}{Lemma}

\newtheorem{definition}{Definition}
\newtheorem{assum}{Assumption}
\newtheorem{remark}{Remark}

\newtheorem{example}{Example}

\newcommand{\I}{{\mathbb{I}}}

\begin{document}
%
\title{\huge{An Improved Distributed Nonlinear Observer for Leader-Following Consensus Via Differential Geometry Approach}}

\author{Haotian Xu,
        Jingcheng~Wang,
        Bohui~Wang,~\IEEEmembership{Member,~IEEE},
        Hongyuan Wang,
        Ibrahim Brahmia
\thanks{This work is supported by National Natural Science Foundation of China (No.61533013, 61633019, 61433002), Key projects from Ministry of Science and Technology (No.2017ZX07207003, 2017ZX07207005), Shaanxi Provincial Key Project (2018ZDXM-GY-168) and Shanghai Project (17DZ1202704). The corresponding author: Jingcheng~Wang.}
\thanks{H. Xu, H. Wang and I. Brahmia are with Department of Automation, Shanghai Jiao Tong University, Key Laboratory of System Control and Information Processing, Ministry of Education of China, Shanghai, 200240~(e-mail: xuhaotian\_1993@126.com; wanghongyuan@sjtu.edu.cn)}
\thanks{J. Wang is with Department of Automation, Shanghai Jiao Tong University, Key Laboratory of System Control and Information Processing, Ministry of Education of China, Shanghai, 200240, also the professor in Autonomous Systems and Intelligent Control International Joint Research Center, Xi'an Technological University, Xi'an, Shaanxi, China, 710021~(e-mail:jcwang@sjtu.edu.cn)}
\thanks{B. Wang is with School of Aerospace Science and Technology, Xidian University, Xi'an 710071, China~(e-mail: wang31aa@126.com)}}

\markboth{Submitted to IEEE Transactions on Systems, Man and Cybernetics: Systems}%
{Shell \MakeLowercase{\textit{et al.}}: Bare Demo of IEEEtran.cls for Journals}
%



\maketitle


\begin{abstract}

This paper is concerned with the leader-following output consensus problem in the framework of distributed nonlinear observers. In stead of certain hypotheses on the leader system, a group of geometric conditions is put forward to develop a novel distributed observer strategy with less conservatism, thereby definitely improving the applicability of the existing results. To be more specific, the improved distributed observer can precisely handle consensus problems for some nonlinear leader systems which are invalid for the traditional strategies with the certain assumption, such as Elastic Shaft Single Linkage Manipulator (ESSLM) systems and most of first-order nonlinear systems.


We prove the exponential stability of our distributed observer by proposing two pioneered lemmas to show the quantitative relationship between the maximum eigenvalues of two matrices appearing in Lyapunov type matrices. Then, a partial feedback linearization method with zero dynamic proposed in differential geometry is employed to design a purely decentralized control law for the affine nonlinear multi-agent system. With this advancement, the existing results can be regarded as a specific case owing to that the followers can be chosen as an arbitrary minimum phase affine smooth nonlinear system. We also prove the certainty equivalence principle for the distributed observer-based control law including novel distributed nonlinear observer and improved purely decentralized control law. Our method is illustrated by ESSLM system and Van der Pol system as leader.
\end{abstract}


%
\IEEEpeerreviewmaketitle

~\\
~\\
\emph{This work has been submitted to the IEEE for possible publication. Copyright may be transferred without notice, after which this version may no longer be accessible.}

\section{Introduction}
%
%
%
%
\IEEEPARstart{M}{ulti-agent} system has been widely studied over the last decades. Not only for the leader-following consensus problem \cite{Wang2017Cooperative,Cai2014The,Zhongkui2010Consensus,Wang2016Joining}, but also for the cooperative output regulation problem\cite{Liu2012Robust,Huang2017The,
BohuiCooperative}. An obstacle in the research of multi-agent system is the communication constraints between followers. It means that a follower may not obtain the information from leader or other followers. One follower can only obtain the information from some specific followers, such as its neighbor. In order to handle the communication constraints, a method named \emph{distributed observer} was proposed \cite{Hong2008Distributed,Cai2017The}.


According to the literatures in recent ten years, distributed observer can be divided into two categories. The first kind of distributed observer does not include the leader, and each of its local observer needs to be able to observe all the states of the whole system by using its own output measurements and the state estimates of its neighbors via communication network \cite{Rego2018A,Kim2016Distributed,Park2017Design,Wang2017A,Han2018Towards,Han2017A}. In the second kind of distributed observer, the follower's local observer is designed to estimate the leader's states. In this scene, only a part of the followers can obtain the actual states of leader, while the other followers may achieve accurately state estimate of leader through the information interaction in the communication graph.

Although both of them are called distributed observer, in fact, the first kind of distributed observer is closer to the problem of distributed filtering. Therefore, though distributed filtering has been studied for decades \cite{Khan2008Distributing,Carli2008Distributed,Xie2012Fully,Zhou2013Coordinated,Talarico2014Distributed,Battistelli2017Distributed,Liu2017Asynchronous,Wu2018A,Zeng2017Distributed}, its research method is not suitable for the second kind of distributed observer. The latter, as the main research object of this paper, has become a research hotspot only in the last decade. \cite{YangObserver} has developed observer-based event-triggered leader-following control for a class of linear multi-agent systems. \cite{HongDistributed} has been concerned with a leader-following problem for a multi-agent system with a switching interconnection topology, and leader-following consensus problem for a class of uncertain nonlinear multi-agent systems with linear leader under jointly connected directed switching networks has been investigated by \cite{LiuAdaptive}. A distributed control law, the so-called distributed observer approach or distributed observer-based framework has been initialed proposed by \cite{Su2012Cooperative}.

\begin{figure}[!t]
  \centering
  \includegraphics[width=8cm]{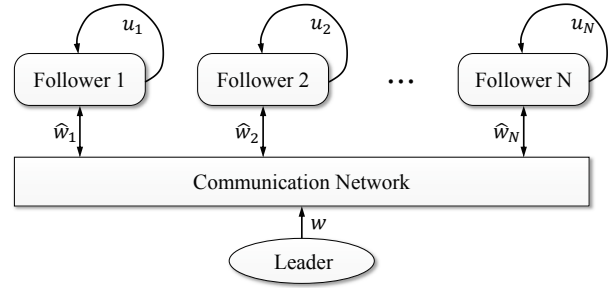}\\
  \caption{The structure of the distributed observer-based control. $\omega$ is the actual state of leader and all the followers need to obtain state estimate $\hat{\omega}_i$ of leader via communication network. The control law $u_i$ of each follower is generated by their own observer information $\hat{\omega}_i$.}\label{structure}
\end{figure}

However, most of these researches focus on linear systems or linear leader system \cite{LiuAdaptive}, only a small part of them take nonlinear system into account, for example, see attitude control of rigid body based on distributed observer \cite{Cai2014A}. Recently, \cite{LiuA} has designed a distributed nonlinear observer for a class nonlinear system, and successfully established the same observer-based distributed control framework for leader-following problem in terms of nonlinear system as that of linear system. This framework, see Figure \ref{structure}, is a general way to design the distributed observer and tracking controller under the leader-following information, which includes a distributed nonlinear observer for leader system and a group of purely decentralized control law assigned to follower systems. The observer-based distributed control law is formed by replacing the leader's actual states which cannot be obtained in purely decentralized control law with the state estimate generated by the local observer.

The most formidable task of this framework is to ensure the stability of the distributed observer, which is also an important basis to ensure that the distributed control law satisfies the certainty equivalency principle \cite{HuangControl}. The stability of \cite{LiuA}'s distributed observer is achieved on basis of Assumptions that the leader system is globally bounded, and ought to meet \emph{Taylor conditions}, which requires that the nonlinear term of leader's Taylor expansion around origin in first-order fashion is in the form of $p(w)w$, where $w$ is state vector and $p(w)$ is a diagonal matrix with all its diagonal entries no more than $0$. However, Taylor conditions for leader system is too strict to be acceptable even by some common systems. For example, see a simple system $\dot{w}=-\sin w$. The diagonal entries $-\sin w$ of its nonlinear term around origin cannot be guaranteed to be less than $0$.

Motivated by this problem, this paper imposes constraints on leader system by a group of geometric conditions (conditions of OCF) \cite{Krener1983Linearization,Isidori2003Nonlinear,Krener1985,Xia2006Nonlinear,Ramdane2016Partial,Noh2004Nonlinear,Sassano2018A} instead of Taylor conditions, and develops a new algorithm to derive the distributed nonlinear observer and tracking controller. Although there are a few papers that design distributed nonlinear observer based on differential geometry method, such as \cite{Silm2019A} and \cite{XuNolinear2019}, they all study the first kind of distributed observer. To the author's knows, this is the first paper to construct the second kind of distributed nonlinear observer with geometry condition. Based on this method, the application range of distributed nonlinear observer with regard to \cite{LiuA}'s distributed control framework can be extended. For example, all of the first order nonlinear smooth system can be accepted by our distributed observer but most of them cannot satisfy \cite{LiuA}'s assumption. Fortunately, some practical systems that fails to meet the Taylor conditions, such as Elastic Shaft Single Linkage Manipulator System (ESSLM), can satisfy our geometric conditions. Moreover, as the first paper in the research of nonlinear distributed control framework, \cite{LiuA} designs a purely decentralized control law for a completely controllable single input follower. Consequently, it is the second purpose of this investigation to study the decentralized control law for more general follower systems, such as multi input systems or incompletely controllable systems.

This investigation, whether it aims to expand the application scope of distributed nonlinear observer or it improves the design method of purely decentralized control law, is by no means trivial. The challenges mainly comes from: 1) The system meeting geometric conditions can be transformed into quasi linearized form via a diffeomorphism; then how to make full use of the property of diffeomorphism to construct an observer with quasi linear error dynamics and deal with the nonlinear term in the quasi linearized observer? 2) During to the introducing of geometric conditions, a relationship between the maximum eigenvalues of two matrices appearing in Lyapunov form matrix is necessary to guarantee the stability of distributed observer; how to find this relationship is one of the key issue in this investigate. 3) How to find a diffeomorphism to make the follower system and the tracking error system between leader and follower have the same zero dynamics? It is the key to solve the leader-following consensus problem when the follower is the minimum phase affine nonlinear system.

The main contributions of this paper consist of the following five aspects:
\begin{itemize}
  \item We constraint nonlinear leader system with geometric conditions instead of Taylor conditions, which enlarges the application range of their framework about nonlinear distributed control law.
  \item A novel distributed nonlinear observer based on differential geometry is proposed, and it is proved to achieve exponential stability for all output bounded affine nonlinear system which meets geometric conditions.
  \item In order to prove the stability of the novel distributed observer, we analysis the relationship between the maximum eigenvalues of two matrices appearing in Lyapunov type matrices carefully with the help of inequality analysis and matrix theory. And describe this relationship in form of inequality.
  \item For certainty systems, the using of differential geometry based novel distributed observer leads to less conservatism. For example, for Van der Pol system, our method can obtain globally convergent distributed observer, rather than converging only in a compact set containing the origin.
  \item This paper develops purely decentralized control law based on zero dynamic theory, which enable the selection range of follower system to be expanded from fully controllable affine nonlinear system to minimum phase affine nonlinear system.
\end{itemize}

This paper is organized as follows. Some notations and mathematical tools throughout this paper are summarized in Section \ref{sec2}. The motivation and improvement of our novel distributed observer are detailed with two examples, including an actual system ESSLM, in \ref{sec3}. Section \ref{sec4} proves the stability of our novel distributed observer with two eigenvalues Lemma. The improved purely decentralized control law and the proof of the principle of definite equivalence are shown in Section \ref{sec5}. In Section \ref{sec6}, the results of our method are simulated with the Van del Pol system and ESSLM system as the leader. Finally, Section \ref{sec7} concludes this paper.

\section{Preliminaries and Problem Formulation}
\label{sec2}

\subsection{Notation}

$I_N\in R^{N\times N}$ denotes an identity matrix. $1_N$ is an $N$ dimensional column vector with all its entries equaling $1$. $A^*$ and $A^T$ denote the conjugate transposition and the transposition of matrix $A$ respectively. For some column vectors $a_i$, denote $col\{a_1,a_2,\cdots,a_N\}$ as a column vector $[a_1^T,a_2^T,\cdots,a_N^T ]^T$. For some matrices $A_i$, $diag\{A_1,A_2,\cdots, A_N\}$ represents a block diagonal matrix. Let $\bar{\sigma}(P)$ be the maximum real of all eigenvalues of $P$, and $\underline{\sigma}(P)$ be the minimum real of all eigenvalues. Particularly, if $P$ Hermite matrix, $\bar{\sigma}(P)$ and $\underline{\sigma}(\cdot)$ represent the maximum and minimum eigenvalue respectively. $\otimes$ denotes Kronecker product with a property $(A\otimes B)(C\otimes D)=AC\otimes BD$. $\|\cdot\|$ is denoted as the 2-norm of matrix or vector.
\subsection{Graph theory}
A directed graph is usually expressed as $\mathcal{G}=(\mathcal{V},\mathcal{E},\mathcal{A}_{\mathcal{G}})$, where $\mathcal{V}$ is the node set including $v_1,v_2,\cdots,v_N$, $\mathcal{E}$ is the arc set, and $\mathcal{A}_{\mathcal{G}}=[a_{ij}]$ is an adjacency matrix of $\mathcal{G}$. Herein, $\mathcal{G}$ is assumed that there are no repeated arcs and no self loops. We denote $a_{ij}=1$ if there is an arc from $v_j$ to $v_i$, denoted as $(v_j,v_i)$, otherwise $a_{ij}=0$. A directed path from node $i$ to node $j$ is a sequence of arcs, expressed $\{(v_i,v_k),(v_k,v_l),\cdots,(v_m.v_j)\}$. The in-degree matrix is defined as $\mathcal{D}_{\mathcal{G}}=diag\{d_i\}$ with $d_i=\sum_{j=1}^Na_{ij}$. The Laplacian matrix of $\mathcal{G}$ is in form of $\mathcal{L}=\mathcal{D}_{\mathcal{G}}-\mathcal{A}_{\mathcal{G}}$. An extended graph is $\bar{\mathcal{G}}=(\bar{\mathcal{V}},\bar{\mathcal{E}},\mathcal{A}_{\mathcal{G}})$, where $\bar{\mathcal{V}}=\mathcal{V}\cup v_0$ with $v_0$ being the node associated with leader, $\bar{\mathcal{E}}$ includes all the arcs in $\mathcal{E}$ and all the arcs between $v_0$ and $\mathcal{E}$. Denote $B=diag\{b_i\}$ where $b_i=1$ if $(v_0,v_i)\in\bar{\mathcal{E}}$, otherwise $b_i=0$. One may know from references \cite{Hu2007Leader,Zhongkui2010Consensus,Qu2009Cooperative} that $\mathcal{L}+B$ is a semi-defined positive matrix if $v_0$ is a globally reachable node. The leader node is a \emph{global reachable} node \cite{Hu2007Leader} if and only if for each $v_i\in\bar{\mathcal{V}}$, there is a directed path from $v_0$ to $v_i$.

\subsection{Differential geometry}

One can refer the knowledge of this subsection and next subsection to \cite{LiDianpu2006} or \cite{Isidori2003Nonlinear}. The Lie derivative of a smooth function $h$ along the vector field $f$ is defined as $L_fh=\frac{\partial h}{\partial x^T}f(x)$. Moreover, if there is a dual vector field $\omega$ belonging to dual tangent space, the Lie derivative of dual vector field $\omega$ along to the vector field $f$ is expressed as $L_f\omega=f^T(\frac{\partial\omega^T}{\partial x^T})^T+\omega\frac{\partial f}{\partial x^T}$. $[f,g]=\frac{\partial g}{\partial x^T}f-\frac{\partial f}{\partial x^T}g$ denotes the Lie bracket between two vector fields. $[f,g]$ is also denoted as $ad_fg$, and $ad_f^kg=[f,ad_f^{k-1}g]$ for $k\geq 2$ is the notation of higher-order Lie bracket, where $ad_f^0g=g$, and $ad_f^1g=ad_fg$. A \emph{distribution} $\mathcal{D}$ is spanned by a group of vector fields $X_1,\cdots,X_d$, i.e., $\mathcal{D}=span\{X_1,X_2,\cdots,X_d\}$. We call $\mathcal{D}$ is involutive if $[X_i,X_j]\in\mathcal{D}$ for $\forall X_i, X_j\in \mathcal{D}$.

\subsection{Zero dynamics}\label{sec2.4}

Given a $n$ dimensional SISO system:
\begin{align}
\dot{x}&=f(x)+g(x)u,\\
y&=h(x).
\end{align}
The \emph{relative degree} of this system is $r$ in a neighborhood $\mathcal{U}$ around a given point $x^0$ if $L_gL_f^kh(x)=0$ for all $x\in\mathcal{U}$ and $0\leq k<r-1$, and $L_gL_f^{r-1}h(x^0)\neq 0$. This system can be transformed (see the coordinate transformation $z=\Phi(x)$ in \cite[Chapter 10]{LiDianpu2006}) into a normal form
\begin{align}
\dot{z}_i&=z_{i+1}, i=1,2,\cdots,r-1,\\
\dot{z}_r&=b(\bar{z},\theta)+a(\bar{z},\theta)u,\\
\dot{\theta}&=\gamma(\bar{z},\theta),\label{internal1}
\end{align}
or a quasi-normal form
\begin{align}
\dot{z}_i&=z_{i+1}, i=1,2,\cdots,r-1,\\
\dot{z}_r&=b(\bar{z},\theta)+a(\bar{z},\theta)u,\\
\dot{\theta}&=\gamma(\bar{z},\theta)+\rho(\bar{z},\theta)u,\label{internal2}
\end{align}
where $\theta=col\{z_{r+1},\cdots,z_{n}\}$, $\bar{z}=col\{z_1,\cdots,z_r\}$, $a(z)=L_gL_f^{r-1}h(\Phi^{-1}(z))$, and $b(z)=L_f^rh(\Phi^{-1}(z))$. Equations (\ref{internal1}) and (\ref{internal2}) are denoted as \emph{internal dynamics} of the original system. By restricting the internal dynamics on \emph{zero dynamics space} (i.e., set $z_1=\cdots=z_r=0$), one may get
\begin{equation}
\dot{\theta}=\gamma(0,\theta), \label{zero1}
\end{equation}
for normal form and may get
\begin{equation}
\dot{\theta}=\gamma(0,\theta)-\rho(0,\theta)\frac{b(0,\theta)}{a(0,\theta)}, \label{zero2}
\end{equation}
for quasi-normal form. (\ref{zero1}) and (\ref{zero2}) are called \emph{zero dynamic} of the original nonlinear system. The relative degree and coordinate transformation about MIMO system will be introduced in Section \ref{sec4}, the corresponding definition of zero dynamics can be defined in a similar way. A nonlinear system is called \emph{minimum phase system} if it has a stability zero dynamic.
\subsection{Problem Formulation}
Consider a nonlinear multi-agent systems with all subsystems being in the form of affine nonlinear as follow
\begin{align}
\dot{x}_i&=f_i(x_i,w)+g_i(x_i)u_i,\label{followsys1}\\
y_i&=h_i(x_i),\ \ i=1,2,\cdots,N, \label{followsys2}
\end{align}
where $x_i\in\mathbb{R}^{n_i}$, $u_i\in\mathbb{R}^m$, $y_i\in\mathbb{R}^r$ are the states, control inputs, and measurement outputs of the $i$th subsystem, or named follower. $f_i(\cdot)$, $g_i(\cdot)$, $h_i(\cdot)$ are smooth vector value functions. And the variable $w\in\mathbb{R}^s$ is generated by an external system, or called leader system, which is an autonomous system,
\begin{align}
\dot{w}&=p(w), \label{leadersys1}\\
y_0&=q(w), \label{leadersys2}
\end{align}
where $p(\cdot)$ and $q(\cdot)$ are smooth vector value nonlinear functions and $y_0\in\mathbb{R}^r$ is the output of the external system. The problem is to design a distributed control law to let the measurement output of subsystem track the output of leader system, i.e.
\begin{equation}
\lim_{t\to\infty}y_0(t)-y_i(t)=0,\ \ \ i=1,2,\cdots,N. \label{object}
\end{equation}

Similar to the problem background of \cite{LiuA}, all followers can get the information of their neighbors' agents via a communication network, but only a part of followers can obtain the real states of leader system. Since the leader-follower problem needs to add the leader's states to the followers' control law, which requires all followers to estimate the leader's state by their own and neighbor's information. Aiming at this problem, a frame of observer-based distributed control for nonlinear system proposed in \cite{LiuA}. This frame contains three aspects. First is to design a distributed observer based on communication graph and show that whether it is existing for the studied nonlinear system. Secondly, one should design purely decentralized control law for every subsystem. Finally, the distributed observer and the purely decentralized control law constitute the distributed control law together.

In the rest of this paper, we will focus on how to improve the distributed nonlinear observer enlarge the application range of \cite{LiuA}'s distributed control framework, and how to design the distributed control law when the follower systems are minimum phase system.
\section{Motivation and Improvement}
\label{sec3}
The important premise of the distributed control framework described in Section \ref{sec2} is the stability of the distributed observer. However, the nonlinear leader system makes researchers have to limit the form of the system, otherwise it is difficult to guarantee the stability of the distributed observer.

In this section, we will first review the assumptions added to the leader system in \cite{LiuA}. Then the geometric conditions will be proposed to constrain the leader system. We will show by two examples that our conditions can enlarge the application range of \cite{LiuA}'s distributed control framework. To move on, a new distributed observer based on geometric conditions will be designed in the Section \ref{sec3.3}, and we will give our main result at the end of this section that the new distributed observer can achieve stability for the output bounded nonlinear system which meet the geometric conditions.

\subsection{Existing results and Motivation}
The $i$th local observer of distributed observer in \cite{LiuA} is introduced as:
\begin{equation}
\dot{\hat{\omega}}_i=p(\hat{\omega})+c\sum_{i=1}^Na_{ij}(\hat{\omega}_j-\hat{\omega}_i).
\end{equation}
where $\hat{\omega}_i$ represents the state estimate of leader system given by the $i$th follower system, and $c$ is the coupling gain. The stability of this distributed observer is based on the following three basic assumptions:

\begin{assum}\label{assum1}
The dynamic of leader system (\ref{leadersys1})(\ref{leadersys2}) is output bounded.
\end{assum}

\begin{assum}\label{assum2}
The leader node in communication network $\bar{\mathcal{G}}$ is assumed to be a globally reachable node. This condition is equivalent to suppose the communication network $\bar{\mathcal{G}}$ has a $v_0$-spanning tree.
\end{assum}

\begin{assum}\label{assum3}
The dynamic function $p(w)$ of leader system is supposed in the following form around the origin, i.e.,
\begin{equation*}
  \left.p(w)=\frac{\partial p}{\partial w}\right|_{w=0}w+p_2(w)w,
\end{equation*}
where $p_2(w)=diag\{d_1(w),d_2(w),\cdots, d_s(w)\}$ with all its diagonal entries being less than zero, i.e., $d_i(w)\leq 0$ for all $i=1,2,\cdots,s$.
\end{assum}

However, Assumption \ref{assum3}, so-called Taylor condition, is too strict to be fulfilled, because few systems can meet this requirement that all the diagonal elements of $p_2(\omega)$ are not positive. For example, the nonlinearity of the following two systems is not very strong, but neither of them can satisfy Assumption \ref{assum3}.

\begin{example}\label{example1}
Consider a numerical example
\begin{equation}
\begin{split}
\dot{w}_1&=-w_1w_2^2+w_3,\\
\dot{w}_2&=-w_1-w_2w_4,\\
\dot{w}_3&=-w_3w_4^2+w_2,\\
\dot{w}_4&=-w_3,\\
y_{01}&=w_2,\ y_{02}=w_4.
\end{split}
\end{equation}
Following Assumption \ref{assum3}, this system can be rewritten as
\begin{equation}
\dot{w}=\begin{pmatrix}0&0&1&0\\-1&0&0&0\\0&1&0&0\\0&0&-1&0\end{pmatrix}w+\begin{pmatrix}-w_2^2&0&0&0\\0&-w_4&0&0\\0&0&-w_4^2&0\\0&0&0&0\end{pmatrix}w.
\end{equation}
Assumption \ref{assum3} cannot be accepted by this system owing to one cannot guarantee that $\omega_4\leq 0$. Hence, the observer proposed in \cite{LiuA} is not able to be applied for this system.
\end{example}
\begin{figure}
  \centering
  \includegraphics[width=6cm]{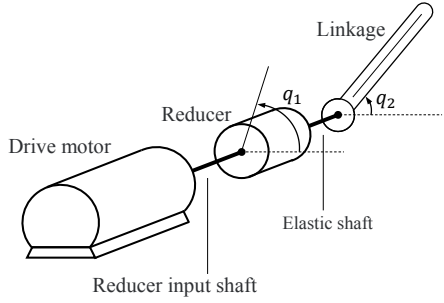}\\
  \caption{Elastic Shaft Single Linkage Manipulator System}\label{linkage}
\end{figure}

\begin{example}\label{example2}
Consider the Elastic Shaft Single Linkage Manipulator System (ESSLM), see Figure \ref{linkage}. Let the length of the Linkage be $2d$ and the mass be $m$. The angular displacement of reducer input shaft and reducer output shaft are $\omega_1$ and $\omega_1/\varpi$ respectively, where $\varpi$ is the transmission ratio of reducer. Denote the angular displacement of the Linkage is $\omega_2$, then the torque at both ends of elastic shaft is $\mathcal{K}(\omega_2-\omega_1/\varpi)$ with $\mathcal{K}$ representing the torsional elastic coefficient. We denote the Viscosity friction coefficient and Rotational inertia of Motor are $F_1$ and $J_1$ respectively and further suppose the Viscosity friction coefficient and the Rotational inertia of Reducer are $F_2$ and $J_2$ respectively. Then the system equation of ESSLM can be introduced as:
\begin{equation}\label{esslm}
\begin{split}
\begin{bmatrix}\dot{\omega}_1\\ \dot{\omega}_2\\ \dot{\omega}_3\\ \dot{\omega}_4\end{bmatrix}
&=p(\omega)\triangleq\begin{bmatrix}\omega_3\\\omega_4\\-\frac{\mathcal{K}}{J_1N^2}\omega_1+\frac{\mathcal{K}}{J_1N}\omega_2-\frac{F_1}{J_1}\omega_3\\
\frac{\mathcal{K}}{J_2N}\omega_1-\frac{\mathcal{K}}{J_2}\omega_2-\frac{mgd}{J_2}\cos{\omega_2}-\frac{F_2}{J_2}\omega_4\end{bmatrix},\\
y&=q(\omega)=\omega_2.
\end{split}
\end{equation}
Note that only $\omega_4$ dynamic in this system is nonlinear. Nevertheless, it still cannot satisfy Assumption \ref{assum3}. Actually, it can be rewritten as
\begin{align}
&\dot{\omega}=\left.\frac{\partial p}{\partial\omega}\right|_{\omega=0}\omega+p_2(\omega)w\notag\\
=&-\begin{bmatrix}0&0&1&0\\0&0&0&1\\ \frac{\mathcal{K}}{J_1\varpi^2}&-\frac{\mathcal{K}}{J_1\varpi}&\frac{F_1}{J_1}&0\\
-\frac{\mathcal{K}}{J_2\varpi}& \frac{\mathcal{K}}{J_2}&0&\frac{F_2}{J_2}\end{bmatrix}\omega+
\begin{bmatrix}0&0&0&0\\0&0&0&0\\0&0&0&0\\0&0&0&d_4(\omega)\end{bmatrix}\omega,
\end{align}
where $d_4(\omega)=-\frac{b_3\cos{\omega_2}}{\omega_4}$ does not fulfill $d_4(\omega)\leq 0$.
\end{example}

\subsection{geometric conditions}\label{sec3.2}
It is known from the previous subsection that too strict system form will limit the application range of distributed control framework. In order to make this framework be widely used in some common nonlinear systems, we introduce the geometric conditions of \emph{observable canonical form} (OCF). Then the distributed observer is designed for the system which can meet the geometric conditions. The observable canonical form can be described as
\begin{align}
\dot{\eta}_0&=A_0\eta_0+a(y_0), \label{linearleader}\\
y_0&=C\eta_0,\notag
\end{align}
where
\begin{gather}
A_0=diag\left\{\begin{bmatrix}0&0\\I_{k_1-1}&0\end{bmatrix},\begin{bmatrix}0&0\\I_{k_2-1}&0\end{bmatrix},\cdots,\begin{bmatrix}0&0\\I_{k_r-1}&0\end{bmatrix}\right\}, \label{A0}\\
C=diag\left\{\begin{bmatrix}0_{1\times(k_1-1)}&1\end{bmatrix},\cdots,\begin{bmatrix}0_{1\times(k_r-1)}&1\end{bmatrix}\right\}, \label{C}
\end{gather}
and $r$-tuples $\{k_1,k_2,\cdots,k_r\}$ is called observable relative degree. Some associated references, such as \cite{Krener1985} and \cite{Xia1989Nonlinear}, prove that the nonlinear system can be transformed into OCF by a diffeomorphism on the premise of satisfying some geometric conditions. Before introducing them, we are supposed to define some codistributions. By rewritting the output function $q(w)$ as $q(w)=\left[q_1(w),q_2(w),\cdots,q_r(w)\right]^T$ and letting $k_i$ be the observable relative degree associated with $q_i(w)$, a group of codistributions defined by \cite{Krener1985} can be introduced as
\begin{align}
&\Delta^{\perp}=span\{dL_f^{l}q_j | 1\leq j\leq r, 0\leq l\leq k_j-1\},\label{delta}\\
&\Delta_i^{\perp}=span\{dL_f^{l}q_j\backslash dL_f^{k_i-1}h_i | 1\leq k\leq r, 0\leq l\leq k_i-1\},\notag\\
&i=1,2,\cdots,r.
\end{align}
The geometric conditions for the existence of diffeomorphism are concluded in the following Lemma \cite{Xia1989Nonlinear}.

\begin{lemma}
\label{lemma1}
The leader system is denoted in (\ref{leadersys1})(\ref{leadersys2}) and its observable relative degree is given by $r$-tuples $\{k_1,k_2,\cdots,k_r\}$. Without loss of generality, we assume $k_1\geq k_2\geq\cdots\geq k_r$ and $\sum_{i=1}^Nk_i=s$. Then there is a diffeomorphism $\eta_0=\Theta(\omega)$ defined on a neighborhood $\mathcal{W}$ around a given point $\omega_0$ which can transform (\ref{leadersys1})(\ref{leadersys2}) into OCF (\ref{linearleader}) if and only if\\
(1). The dimension of codistribution $\Delta^{\perp}$ is $n$.\\
(2). $dim\{\Delta_i^{\perp}\}=dim\{\Delta^{\perp}\cap\Delta_i^{\perp}\}$.\\
(3). For the given linear equations:
\begin{align}
    &\left\langle dL_f^{l-1}h_i,\tau_j\right\rangle=\delta_{i,j}\cdot\delta_{l,r_i},\  l=1,\cdots,r_i,\ \ \emph{if}\ i\leq j, \label{le1}\\
    &\left\langle dL_f^{l-1}h_i,\tau_j\right\rangle=\delta_{i,j}\cdot\delta_{l,r_i},\  l=1,\cdots,r_j,\ \ \emph{if}\ i>j,\label{le2}
\end{align}
there exists a group of vector fields $\tau_1,\tau_2,\cdots,\tau_r$ solved by (\ref{le1})(\ref{le2}) s.t.
the communication conditions $[ad_f^l\tau_i,ad_f^k\tau_j]=0$ are satisfied for all $1\leq i,j\leq r$, $0\leq l\leq k_i-1$, and $0\leq k\leq k_j-1$, where $\delta_{i,j}$ is Kronecker delta.
\end{lemma}

In Lemma \ref{lemma1}, conditions (1)-(3) are called geometric conditions. The solving procedure with respect to calculating $\eta_0=\Theta(\omega)$ given by \cite{Xia1989Nonlinear} are demonstrated in Appendix.

\begin{remark}
Condition (2) of Lemma \ref{lemma1} is hard to understand. It can be stated as $dim\{\Delta^{\perp}\cap\Delta_i^{\perp}\}=ik_i+k_{i+1}+\cdots+k_r-1$. See example 1 of \cite{Xia1989Nonlinear} to learn the situation when Condition (2) is not fulfilled. Note that Condition (2) is satisfied if $k_1=k_2=\cdots=k_r$, see in \cite{Xia1989Nonlinear}\cite{Saadi2016Multi}.
\end{remark}

It is easy to verify that the geometric conditions are fulfilled for all of the first-order nonlinear smooth system, such as $\dot{\omega}=-\sin{\omega}$, but most of them cannot satisfy Taylor conditions. Moreover, geometric conditions can also be applied to some high-order nonlinear systems which fails to meet Assumption \ref{assum3}. As a comparison, we will verify that the two examples in the previous subsection meet the geometric conditions.

\begin{example}
Consider the system in Example (\ref{example1}) again. One can check that the observable relative degree satisfies $k_1=k_2=2$, so condition (1)(2) in Lemma \ref{lemma1} are fulfilled. By calculating $\tau_1=[-1,0,0,0]^T$, $ad_p\tau_1=[w_2^2,1,0,0]^T$, $\tau_2=[0,0,-1,0]^T$ and $ad_p\tau_2[1,0,x_4^2,1]^T$, we can verify condition (3) is also satisfied. Thus system (\ref{example1}) can be transformed into observable canonical form.
\end{example}
\begin{example}
Consider ESSLM system again. Herein, we will verify whether it can be transformed into OCF by diffeomorphism. It can be calculated directly that
\begin{equation}
\begin{bmatrix}dq(\omega)\\dL_pq(\omega)\\dL^2_pq(\omega)\\dL^3_pq(\omega)\end{bmatrix}=
\begin{bmatrix}0&1&0&0\\0&0&0&1\\ \frac{\mathcal{K}}{J_2\varpi}& \Upsilon_1(\omega)&0&-\frac{F_2}{J_2}\\
-\frac{\mathcal{K}F_2}{J_2^2\varpi}& \Upsilon_2(\omega)&\frac{\mathcal{K}}{J_2\varpi}&\Upsilon_3(\omega)\end{bmatrix},
\end{equation}
where
\begin{align}
\Upsilon_1(\omega)&=-\frac{\mathcal{K}}{J_2}+\frac{mgd}{J_2}\sin{\omega_2},\notag\\
\Upsilon_2(\omega)&=\frac{\mathcal{K}F_2}{J_2^2}+\frac{mgd}{J_2}\omega_4\cos{\omega_2}-\frac{mgdF_2}{J_2^2}\sin{\omega_2},\notag\\
\Upsilon_3(\omega)&=\frac{F_2^2}{J_2^2}-\frac{\mathcal{K}}{J_2}+\frac{mgd}{J_2}\sin{\omega_2}.\notag
\end{align}
Thus the solution $\tau$ of linear equations (\ref{le1})(\ref{le2}) can be described as $\tau=[0,0,J_2\varpi/\mathcal{K},0]^T$. Then we can further calculate
\begin{align*}
ad_p\tau&=[-J_2\varpi/\mathcal{K},0,F_1J_2\varpi/J_1\mathcal{K},0]^T,\\ ad_p^2\tau&=[-F_1J_2\varpi/J_1\mathcal{K},0,-J_2/J_1\varpi+F_1^2J_2\varpi/J_1^2\mathcal{K},1]^T, \\
ad_p^3\tau&=\left[\frac{J_2}{J_1\varpi}-\frac{F_1^2J_2\varpi}{J_1^2\mathcal{K}},-1,\Gamma,\frac{F_1}{J_1}+\frac{F_2}{J_2}\right]^T,
\end{align*}
where $\Gamma=-\frac{2F_1J_2}{J_1\varpi}+\frac{F_1^3J_2\varpi}{J_1^3\mathcal{K}}$. Since all $ad_p^i\tau, i=0,1,2,3$ are constant vector fields, we can obtain $[ad_p^i\tau,ad_p^j\tau]=0$ for all $i,j=0,1,2,3$.
\end{example}

\begin{remark}
Although not all nonlinear systems satisfying Taylor conditions can satisfy geometric conditions, these two examples also show that there are a large number of systems that can satisfy geometric conditions but can not satisfy Taylor conditions. Furthermore, in Section \ref{sec4}, we will prove that the observer designed based on geometric conditions can converge exponentially at any speed. In other words, our distributed observer is not inferior to \cite{LiuA}'s observer in performance, and it can also be applied to some nonlinear systems that \cite{LiuA}'s observer cannot be competent for. Hence, the application range of the distributed control frame proposed by \cite{LiuA} can be extended by this paper.
\end{remark}

\subsection{New Distributed observer}
\label{sec3.3}
Suppose the leader system can meet the geometric conditions list in Lemma (\ref{lemma1}). Then our new distributed observer can be designed on the basis of OCF (\ref{linearleader}). In this scene, the $i$th local observer to leader system can be introduced in the form:
\begin{gather}
\dot{\hat{\eta}}_i=A_0\hat{\eta}_i+a(C\hat{\eta}_i)+cF\varsigma_i, \label{dov1}\\
\varsigma_i=\sum_{i=1}^Na_{ij}\left(\hat{\eta}_j-\hat{\eta}_i\right)+b_i\left(\hat{\eta}_i-\eta_0\right),\label{dov2}
\end{gather}
where $F$ is the LQR gain matrix, $c$ is the coupling gain, and $\varsigma_i$ is the global error dynamic.

The observer error of $i$th subsystem is defined as $e_i=\hat{\eta}_i-\eta_0$, which includes a quasi-linear error dynamic:
\begin{equation}
\dot{e}_i=A_0e_i+a(\hat{y}_i)-a(y_0)+cF\varsigma_i.
\end{equation}
By setting $e=col\left\{e_1,e_2,\cdots,e_N\right\}$, $\tilde{a}_i=a(\hat{y}_i)-a(y_0)$, and $\tilde{a}=col\left\{\tilde{a}_1,\tilde{a}_2,\cdots,\tilde{a}_N\right\}$, we can rewrite the error dynamic in a compact form:
\begin{align}
\dot{e}&=\left(I_N\otimes A_0\right)e-c\left(I_N\otimes F\right)\left(\left(\mathcal{L}+B\right)\otimes I_N\right)e+\tilde{a} \notag \\
&=\left(I_N\otimes A_0-c\left(\mathcal{L}+B\right)\otimes F\right)e+\tilde{a}. \label{error}
\end{align}
Denote $\mathcal{M}=I_N\otimes A_0-c\left(\mathcal{L}+B\right)\otimes F$. An intuitive fact could be noticed that the properties of $\mathcal{M}$ have an important influence on the stability of the error system (\ref{error}). Fortunately, some results \cite{Zhang2011Optimal} in multi-agent problem can be found to help us understand the properties of $\mathcal{M}$ and design LQR gain matrix $F$.
\begin{lemma}
\label{lemma2}
Suppose $\lambda_i, i=1,2,\cdots,N$ are the eigenvalues of $\mathcal{L}+B$, then matrix $\mathcal{M}$ is Hurwitz if and only if $A_0-c\lambda_iF$ is Hurwitz for all $i$.
\end{lemma}

This lemma indicates that the stability of $\mathcal{M}$ depends on the structure of communication graph $\mathcal{G}$. Hence, one is supposed to factor the effect of $\lambda_i$ into the mix of designing $F$. Lemma \ref{lemma3} \cite{Zhang2011Optimal} chooses the gain matrix $F$ based on the LQR optimal control. One may refer the proof of lemma \ref{lemma3} in \cite{Zhang2011Optimal}.

\begin{lemma}
\label{lemma3}
Suppose $Q,R$ are symmetric positive definite matrices, and choose $F$ from
\begin{equation}
F=P_1R^{-1}, \label{F}
\end{equation}
where the symmetric positive definite matrices $P_1$ is solved by algebraic Riccati equation
\begin{equation}
AP_1+P_1A^T+Q-P_1R^{-1}P_1=0.
\end{equation}
Then $\mathcal{M}$ is Hurwitz if the coupling gain $c$ satisfies
\begin{equation}
c\geq\frac{1}{2\underline{\sigma}(\mathcal{L}+B)}.\label{c1}
\end{equation}
\end{lemma}

Now we can give one of the main results (Theorem \ref{thm1}) of this paper. This theorem guarantees that the distributed control frame for nonlinear leader-following consensus proposed by \cite{LiuA} can be applied to a class of output bounded nonlinear leader systems which satisfies geometric conditions.

\begin{theorem}\label{thm1}
Suppose the nonlinear leader system (\ref{leadersys1})(\ref{leadersys2}) satisfies Assumptions \ref{assum1} and \ref{assum2} and the pair $(p(\omega),q(\omega))$ meets the geometric conditions proposed in Lemma \ref{lemma1}. Then there exists a coupling gain $c$ satisfying (\ref{c1}) such that the state estimate generated by distributed observer (\ref{dov1})(\ref{dov2}) converges exponentially to actual state of the leader system at arbitrary speed.
\end{theorem}

\section{Proof of Theorem \ref{thm1}}\label{sec4}
In order to deduce the conclusion of Theorem \ref{thm1} in a more accurate way, we are going to prove two Lemmas in \ref{sec4.1} to reveal the quantitative relationship between the maximum eigenvalues of the two matrices appearing in Lyapunov-form matrix rather than the qualitative relationship given in previous literature \cite{Bidram2014Synchronization}. The main body of the proof of Theorem \ref{thm1} will be given in the \ref{sec4.2} via making full use of the diffeomorphism property.
\subsection{Preparation}\label{sec4.1}
\begin{lemma}
\label{prop1}
Set $A\in\mathbb{R}^{n\times n}$ be an arbitrary matrix, $P\in\mathbb{R}^{n\times n}$ be an symmetric positive definite matrix. Then matrix $T=PA+A^TP$ satisfies
\begin{equation}\label{eigen1}
\bar{\sigma}(T)\leq \sqrt{\bar{\sigma}(A^TA)}\bar{\sigma}(P).
\end{equation}
\end{lemma}

\begin{IEEEproof}
Suppose that $\eta$ is the eigenvector of $T$ corresponding to $\bar{\sigma}(T)$. By letting $\xi=A\eta$, we have
\begin{align*}
&\eta^T(PA+A^TP)\eta=\bar{\sigma}(T)\|\eta\|^2=2\eta^TPA\eta=2\eta^TP\xi.
\end{align*}
According to Cauchy Schwartz inequality and C-F inequality, we can further obtain
\begin{align*}
2\eta^TP\xi &\leq 2\left(\eta^T P\eta\cdot\xi^T P\xi\right)^{\frac{1}{2}}\leq 2\bar{\sigma}(P)\|\eta\|\|\xi\|\\
&=2\bar{\sigma}(P)\|\eta\|\|A\eta\|\leq 2\bar{\sigma}(P)\|\eta\|\|A\|\|\eta\|\\
&=2\|A\|\|\bar{\sigma}(P)\|\eta\|^2=\sqrt{\bar{\sigma}(A^TA)}\bar{\sigma}(P)\|\eta\|^2.
\end{align*}
Hence, we have $\bar{\sigma}(T)\leq \sqrt{\bar{\sigma}(A^TA)}\bar{\sigma}(P)$.
\end{IEEEproof}

\begin{lemma}
\label{prop2}
Suppose $\mathcal{M}\in\mathbb{R}^n$ is a Hurwitz matrix. For a fixed constant $\mu>0$, a unique symmetric positive definite matrix $P$ solved by
\begin{equation}
P\mathcal{M}+\mathcal{M}^TP=-2\mu I_n \label{Lyp1}
\end{equation}
satisfying
\begin{align}
&\bar{\sigma}(P)\bar{\sigma}(\mathcal{M})\leq-\mu. \label{eigen3-1}\\
&\bar{\sigma}(P)\bar{\sigma}\left(\mathcal{M}+\mathcal{M}^*\right)\geq-2\mu.\label{eigen3-2}
\end{align}
Especially, $\bar{\sigma}(P)\bar{\sigma}(\mathcal{M})=-\mu$ if $\mathcal{M}$ is a Hermite matrix.
\end{lemma}

\begin{IEEEproof}
Denote $\lambda_2$ as an eigenvalue of $\mathcal{M}$ with $Re(\lambda_2)=\bar{\sigma}(\mathcal{M})$ and treat $\eta$ as the eigenvector of $\mathcal{M}$ corresponding to $\lambda_2$. By pre-multiplying $\eta^*$ and post-multiplying $\eta$ on (\ref{Lyp1}), we have
\begin{equation*}
\eta^*\left(P\mathcal{M}+\mathcal{M}^*P\right)\eta=-2\mu\eta^*\eta.
\end{equation*}
A natural step can be obtained as:
\begin{equation*}
\lambda_2\eta^*P\eta+\lambda_2^*\eta^*P\eta=2Re(\lambda_2)\eta^*P\eta=-2\mu\eta^*\eta.
\end{equation*}
Note that $\bar{\sigma}(\mathcal{M})<0$ since $\mathcal{M}$ is a Hurwitz matrix. Then by denoting $\lambda_1=\bar{\sigma}(P)$, we can deduce with C-F inequation
\begin{equation*}
-2\mu\eta^*\eta=2\bar{\sigma}(\mathcal{M})\eta^*P\eta\geq 2\bar{\sigma}(\mathcal{M})\bar{\sigma}(P)\eta^*\eta.
\end{equation*}
Consequently,
\begin{equation*}
\bar{\sigma}(P)\bar{\sigma}(\mathcal{M})\leq-\mu.
\end{equation*}
Furthermore, there exists a orthogonal matrix $U$ such that $P=U^T\Lambda U$ because $P$ is a real symmetric matrix, where $\Lambda$ is a diagonal matrix with all eigenvalues of $P$ on its diagonal. Then we calculate by setting $\bar{\mathcal{M}}=U\mathcal{M}U^T$ that,
\begin{align}
&P\mathcal{M}+\mathcal{M}^*P\notag\\
=&U^T\Lambda U\mathcal{M}U^TU+U^TU\mathcal{M}^*U^T\Lambda U\notag\\
=&U^T\left(\bar{\mathcal{M}}^*\Lambda+\Lambda\bar{\mathcal{M}}\right)U\leq \bar{\sigma}(P)U^T
\left(\bar{\mathcal{M}}+\bar{\mathcal{M}}^*\right)U\notag\\
=&\bar{\sigma}(p)\left(\mathcal{M}+\mathcal{M}^*\right).
\end{align}
It yields
\begin{align}
\bar{\sigma}(P)\bar{\sigma}\left(\mathcal{M}+\mathcal{M}^*\right)\geq-2\mu.
\end{align}

In particular, by supposing $\mathcal{M}$ is a Hermite matrix, i.e., $\mathcal{M}=\mathcal{M}^*$, we thus get
\begin{equation*}
\bar{\sigma}(P)\bar{\sigma}\left(\mathcal{M}+\mathcal{M}^*\right)=2\bar{\sigma}(P)\bar{\sigma}(\mathcal{M})\geq-2\mu.
\end{equation*}
Combining this equation and equation (\ref{eigen3-1}), $\bar{\sigma}(P)\bar{\sigma}(\mathcal{M})=-\mu$ can be proved.
\end{IEEEproof}

\begin{remark}
Figure \ref{p1} shows the verification of Lemma \ref{prop1}: all the blue dots are located under red line. Figure \ref{p2} and  Figure \ref{p23} show the conclusions of Lemma \ref{prop2}, where $\mathcal{M}$ is chosen as a random Hurwitz matrix and $P$ is a positive defined matrix solved by (\ref{Lyp1}) with $\mu=2$. We construct one thousand $\mathcal{M}$ and calculate the corresponding $P$. The blue dash line $\bar{\sigma}(P)\bar{\sigma}(\mathcal{M})$ (blue dash line) and red solid line $y=-2$ are plotted in Figure \ref{p2}, and the value of $\bar{\sigma}(P)\bar{\sigma}\left(\mathcal{M}+\mathcal{M}^*\right)$ (blue dash line) and line $y=-4$ (red solid line) are demonstrated in Figure \ref{p23}. These figures show the correctness of equations (\ref{eigen3-1}) and (\ref{eigen3-2}).
\begin{figure}[!t]
  \centering
  \includegraphics[width=8cm]{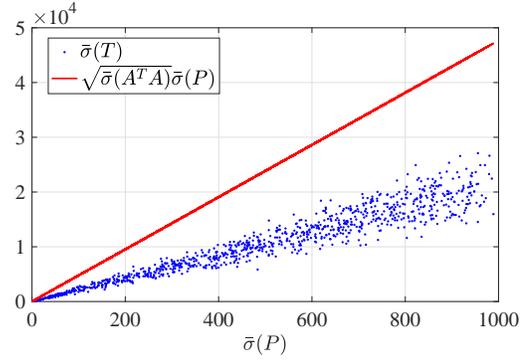}\\
  \caption{Relationship between $\bar{\sigma}(T)$ and $\bar{\sigma}(P)$}\label{p1}
\end{figure}
\begin{figure}[!t]
  \centering
  \includegraphics[width=8cm]{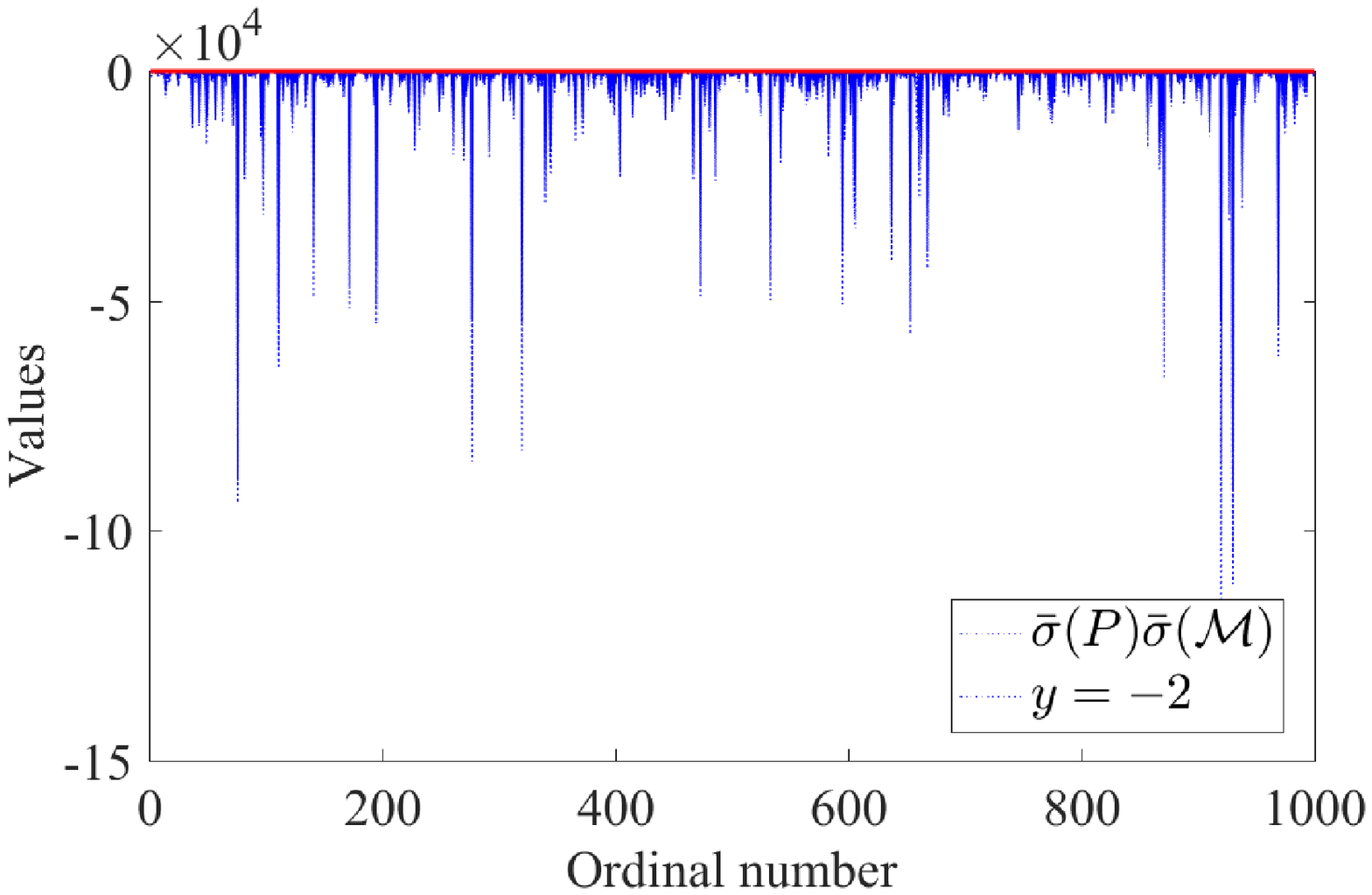}\\
  \caption{Relationship between $\bar{\sigma}(\mathcal{M})$ and $\bar{\sigma}(P)$}\label{p2}
\end{figure}
\begin{figure}[!t]
  \centering
  \includegraphics[width=8cm]{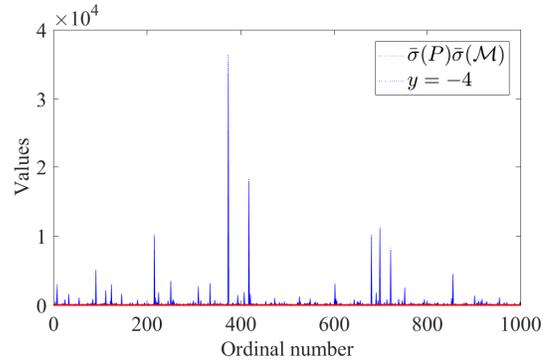}\\
  \caption{Relationship between $\bar{\sigma}(\mathcal{M}+\mathcal{M}^*)$ and $\bar{\sigma}(P)$}\label{p23}
\end{figure}
\end{remark}

\subsection{Main body of the proof}\label{sec4.2}
Before the proof, we are supposed to introduce a definition of \emph{Decreasing Function in Trend}.
\begin{definition}[Decreasing Function in Trend]
A real function $f(x)$ defined on real number field is a decreasing function in trend if for $\forall x_1\in\mathbb{R}$, there exists a $x_2>x_1$, such that $f(x_2)<f(x_1)$.
\end{definition}
Now we prove the conclusion in Theorem \ref{thm1}.
\begin{IEEEproof}
Choose $F$ by the statement of Lemma \ref{lemma3}. Then one can conclude $\mathcal{M}$ is Hurwitz matrix. Thus, there exists a symmetric positive definite matrices $P_2$ such that
\begin{equation}\label{Lyp2}
\mathcal{M}^TP_2+P_2\mathcal{M}=-2\mu I_{sN}.
\end{equation}

Sequentially, Lyapunov function can be chosen as $V(e)=e^TP_2e$, and then we calculate the derivative of $V(e)$ along to error dynamics (\ref{error})
\begin{equation}\label{Vdot}
\dot{V}(e)=e^T\left(\mathcal{M}^TP_2+P_2\mathcal{M}\right)e+2e^TP_2\tilde{a}.
\end{equation}

We know from the steps of observer linearization in appendix that the nonlinear compensation term $a(y_0)$ is the solution of partial differential equations
\begin{equation}\label{pde}
\frac{\partial a(y_0)}{\partial y_{0i}}=b_i(y_0),
\end{equation}
where $b_i(y_0)$ is defined in appendix. However, we can deduce, according to the PDEs in (\ref{pde}) only containing one partial derivative, that every component $a_j(\cdot)$ of $a(\cdot)$ satisfies the differential mean value theorem for each variable $y_{0i}$. It means, for every $a_j(\cdot)$ of $a(\cdot)$ and every $y_{0i}$ of $y_0$, that
\begin{equation}
\tilde{a}_{kj}\left(y_0,\hat{y}_k\right)\leq\sum_{i=1}^r\bar{l}_{ij}\left(\hat{y}_{ki}-y_{0i}\right)=\bar{l}_{j}^T\left(\hat{y}_k-y_0\right),
\end{equation}
where $\bar{l}_{ij}$ is the upper bound of $\partial a_j/\partial y_{0i}$ since output function $y_0$ is bounded. and $\bar{l}_j=\left(\bar{l}_{j1},\bar{l}_{j2},\cdots,\bar{l}_{jr}\right)^T$. Hence, the error of the nonlinear compensation term $\tilde{a}_k\left(y_0,\hat{y}_k\right)$ of $k$th subsystem can be described in a compact equation
\begin{align}
\tilde{a}_k\left(y_0,\hat{y}_k\right)&=\begin{pmatrix}\tilde{a}_{k1}\left(y_0,\hat{y}_k\right)\\\tilde{a}_{k2}\left(y_0,\hat{y}_k\right)\\\vdots\\\tilde{a}_{ks}\left(y_0,\hat{y}_k\right)\end{pmatrix}
\leq\begin{pmatrix}\bar{l}_1^T\\\bar{l}_2^T\\\vdots\\\bar{l}_s^T\end{pmatrix}\left(\hat{y}_k-y_0\right)\notag\\
&\triangleq L\left(\hat{y}_k-y_0\right)=LC\left(\hat{\eta}_k-\eta_0\right).
\end{align}
Therefore, we have
\begin{equation}\label{tilde_a}
\tilde{a}=col\left\{\tilde{a}_1,\tilde{a}_2,\cdots,\tilde{a}_N\right\}\leq \left(I_N\otimes LC\right)e.
\end{equation}
Substituting (\ref{Lyp2}) and (\ref{tilde_a}) into (\ref{Vdot}), we get
\begin{align}\label{Vdot2}
\dot{V}(t)&=e^T\left(-2\mu I_{sN}\right)e+e^TP_2\left(I_N\otimes LC\right)e\notag\\
&=e^T\left(-2\mu I_{sN}\right)e\notag\\
&+\frac{1}{2}e^T\left(P_2\left(I_N\otimes LC\right)+\left(I_N\otimes C^TL^T\right)P_2\right)e\notag\\
&\leq -2\mu \|e\|^2+\frac{1}{2}\kappa\|e\|^2,
\end{align}
where $\kappa=\bar{\sigma}\left(P_2\left(I_N\otimes LC\right)+\left(I_N\otimes C^TL^T\right)P_2\right)$. Moreover, by using of Lemma \ref{prop1}, we know $\kappa\leq\sqrt{\bar{\sigma}(C^TL^TLC)}\bar{\sigma}(P_2)\triangleq\alpha\bar{\sigma}(P_2)$.

Reference \cite{Bidram2014Synchronization} illustrates that one can improve the stability of $\mathcal{M}$ by increasing the coupling gain $c$: Choose a nonsingular matrix $\mathcal{T}$ such that $\mathcal{T}^{-1}(\mathcal{L}+B)\mathcal{T}$ is upper triangular with the eigenvalues ($\lambda_1,\cdots,\lambda_N$) of $\mathcal{L}+B$ on its diagonal. Then $(\mathcal{T}^{-1}\otimes I_n)\mathcal{M}(I_n\otimes\mathcal{T})$ is transformed to $diag\{A_0-c\lambda_iF, i=1,2,\cdots,N\}$. Hence, $\bar{\sigma}(\mathcal{M})\triangleq f(c)$ will decrease in trend with $c$ goes to infinity. Furthermore, $\lim_{c\to\infty}f(c)=-\infty$. By denoting $\mathcal{M}_*=c(\mathcal{L}+B)\otimes F$, we have $\lim_{c\to\infty}\frac{1}{c}(A_0-c\mathcal{M}_*)=-\mathcal{M}_*$, which indicates $\mathcal{M}$ will tends to be $c\left(\mathcal{L}+B\right)\otimes F$ when $c$ tends to infinity.

Now we prove that $\bar{\sigma}(P_2)\bar{\sigma}(\mathcal{M})$ dose not change with $c$ when $c$ is large enough. Actually, we can suppose $c$ is large enough so $\mathcal{M}$ is influenced by $c\left(\mathcal{L}+B\right)\otimes F$ only. Note that $\mathcal{M}$ and $P_2$ are the matrix functions of $c$ denoted as $\mathcal{M}(c)$ and $P_2(c)$ respectively. Then we denote $\mathcal{M}=\mathcal{M}(c)$, $\mathcal{M}'=\mathcal{M}(c+\Delta c)$, $P_2=P_2(c)$ and $P_2'=P_2'(c+\Delta c)$, where $\Delta c$ represents the variation of coupling gain. Thus $\mathcal{M}'$ can be approximately expressed as $c'\mathcal{M}$ when $c$ change to $c+\Delta c$, where $c'=(c+\Delta c)/c$. As a result, $P_2'=P_2/c'$ can be solved by (\ref{Lyp2}). It indicates $\bar{\sigma}(P_2')\bar{\sigma}(\mathcal{M}')=\frac{1}{c'}\bar{\sigma}(P_2)c'\bar{\sigma}(\mathcal{M})=\bar{\sigma}(P_2)\bar{\sigma}(\mathcal{M})$. Hence, $\lim_{c\to\infty}\bar{\sigma}(P_2)\bar{\sigma}(\mathcal{M})=c_1$, where $c_1$ is a constant.

According to Lemma \ref{prop2}, we know $\bar{\sigma}(P_2)\bar{\sigma}(\mathcal{M})\leq -\mu$. Since $\bar{\sigma}(\mathcal{M})$ and $\bar{\sigma}(P_2)$ depend on $c$, we can choose a function $\beta(c)>0$ such that $\bar{\sigma}(P_2)\bar{\sigma}(\mathcal{M})=-\mu-\beta(c)$. Therefore, we can obtain $\lim_{c\to\infty}\beta(c)=c_2$ and $c_2$ is a constant defined by $c_2=-\mu-c_1$. Then the limitation of $\kappa$ can be calculated as
\begin{align}\label{kappa}
&\lim_{c\to\infty}\kappa \leq\lim_{c\to\infty}\alpha\bar{\sigma}(P_2)\notag\\
=&-\alpha\lim_{c\to\infty}\frac{\mu+\beta(c)}{\bar{\sigma}(\mathcal{M})}
=-\lim_{c\to\infty}\frac{\alpha(\mu+\beta(c))}{f(c)}=0.
\end{align}

Therefore, there exists a constant $c^*>0$ such that $\kappa< 4\mu$ for $\forall c>c^*$. It is equivalent to $\dot{V}< 0$. Combining with (\ref{c1}), we know the coupling gain $c$ should satisfy
\begin{equation}\label{c2}
c>\max\left\{\frac{1}{2\underline{\sigma}(\mathcal{L}+B)},c^*\right\}.
\end{equation}
Moreover, for a given $c_0>c^*$, equation (\ref{Vdot2}) can be rewritten as
\begin{align}
\dot{V}(t)&\leq\left(-2\mu+\frac{1}{2}\kappa\right)\|e\|^2\leq \frac{-2\mu+\frac{1}{2}\kappa}{\bar{\sigma}\left(P_2(c_0)\right)}V(t).\notag
\end{align}
So
\begin{align}
V(t)\leq & \exp\left\{\frac{-2\mu+\frac{1}{2}\kappa}{\bar{\sigma}\left(P_2(c_0)\right)} t\right\}V(0)\notag\\
=&exp\left\{\left(\frac{-2\mu}{\bar{\sigma}\left(P_2(c_0)\right)}+\frac{1}{2}\alpha\right)t\right\}V(0).
\end{align}
Since $\lim_{c_0\to\infty}\bar{\sigma}\left(P_2(c_0)\right)=0$, we have
\begin{equation}
\lim_{c_0\to\infty}\frac{-2\mu}{\bar{\sigma}\left(P_2(c_0)\right)}+\frac{1}{2}\alpha=-\infty.
\end{equation}
Thus the error dynamic of this distributed observer can exponential converge to zero at arbitrary speed.
\end{IEEEproof}


\section{Distributed control law for minimum phase affine nonlinear system}
\label{sec5}

In leader-following consensus problem, we only need to control the output related states, the follower system thus need not to be completely controllable. Specifically, the selection range of follower system is expanded from the original completely controllable affine nonlinear system to the minimum phase affine nonlinear system. In this section, we will introduce in detail how to design a purely decentralized control law for the minimum phase follower, especially how to find a differential homeomorphism to make the tracking error system and the follower system have the same zero dynamics. We first introduce the case that the follower system is SISO system, and then extend the problem to the case of MIMO system.

\subsection{Distributed control for SISO system}

Consider an output-tracking problem of leader-following multi-agent system. Leader and follower systems are still in form of (\ref{leadersys1}) and (\ref{leadersys2}) respectively. In this subsection, the leader system is assumed as a single output system, and all followers are derived by SISO nonlinear affine system, i.e., for $i=1,2,\cdots,N$, $y_i$, $y_0$ and the control input signal $u_i$ belong to $\mathbb{R}^1$. In order to study the tracking problem when the follower system is not completely controllable, zero dynamic theory and partial feedback linearization method in differential geometry are employed \cite{Isidori2003Nonlinear}. Within this idea, we propose a purely decentralized control law in which the output of an incompletely controllable follower can track the output of a leader.
\begin{theorem}\label{thm2}
For the $i$th follower system, we assume that it is a minimum phase system and has relative order $r_i$ at $\forall x_i\in\mathbb{R}^{n_i}$. Then there is a coordinate transformation (diffeomorphism) $(\xi_i^T,\theta^T)^T=\Phi_i(x_i)$ such that the tracking error dynamic between the $i$th follower and leader can be described as:
\begin{gather}
\dot{w}=p(w),\label{thm21}\\
\dot{\xi}_i=A_i\xi_i+B_iv_i,\label{thm22}\\
\dot{\theta}_i=\gamma_i\left(\zeta_0+\xi_i, \theta_i\right), \label{thm23}
\end{gather}
where
\begin{equation*}
  A_i=\begin{bmatrix}0&I_{r_i-1}\\0&0\end{bmatrix}\in\mathbb{R}^{r_i\times r_i},\ \ \ B_i=\begin{bmatrix}0&\cdots&0&1\end{bmatrix}\in\mathbb{R}^{r_i},
\end{equation*}
$\zeta_0=col\{q(w),L_pq(w),\cdots,L_p^{r_i-1}q(w)\}$, $v_i$ is a variable named auxiliary control variable, and $\theta_i$ is the internal dynamic of $i$th follower. Furthermore, the leader-following tracking problem (\ref{object}) can be achieved by employing a linear feedback control law
\begin{equation}\label{control1v}
v_i=K_i\xi_i,
\end{equation}
where $K_i$ is a matrix such that $A_i+B_iK_i$ is Hurwitz.
\end{theorem}

\begin{IEEEproof}
Since the $i$th subsystem is a SISO system, the distribution spanned by $g_i(x_i)$ is involutive. Then there is a diffeomorphism $z_i=\Psi_i(x_i)$ \cite{Isidori2003Nonlinear} such that the subsystem can be transformed in normal form
\begin{equation}\label{normal}
\begin{split}
&\dot{z}_{i1}=z_{i2},\\
&\dot{z}_{i2}=z_{i3},\\
&\ \vdots\\
&\dot{z}_{ir_i-1}=z_{it_i},\\
&\dot{z}_{ir_i}=L_{f_i}^{r_i}h_i+L_{g_i}L_{f_i}^{r_i-1}h_iu_i,\\
&\dot{\theta}_i=\gamma_i\left(\zeta_i,\theta_i\right),
\end{split}
\end{equation}
where $\zeta_i=col\{z_{i1},z_{i2},\cdots,z_{ir_i}\}$.

Note that $\Psi=col\{\psi_{i,1},\cdots,\psi_{i,r_i},\psi_{i,r_i+1},\cdots,\psi_{i,n_i}\}$ is constructed by setting $z_{ij}=\psi_{i,j}(x_i)=L_{f_i}^{j-1}h_i$ for $1\leq j\leq r_i$, and choosing $\theta_{i,j}=\psi_{i,j}(x_i), r_i< j\leq n_i$ such that $\Psi_{i*}$ is nonsingular and $L_{g_i}\psi_{i,j}(x_i)=0, r_i< j\leq n_i$. Hence, for $\forall x_i\in\mathbb{R}$, we yield
\begin{align}
n_i&=rank\{\Psi_{i*}\}\notag\\
&=rank\{dh_i,\cdots,dL_{f_i}^{r_i-1}h_i,d\psi_{i,r_i+1},\cdots,d\psi_{i,n_i}\}.
\end{align}

Denote $\varepsilon^{(k)}_i(t)$ as the $k$-order derivative of $\varepsilon_i(t)$ with $\varepsilon_i(t)=y_i(t)-y_0(t)$. Then we can derive from the definition of relative degree that
\begin{gather}
\varepsilon_i^{(1)}=L_{f_i}h_i-L_pq,\label{epsilon1}\\
\varepsilon_i^{(2)}=L_{f_i}^2h_i-L_p^2q,\label{epsilon2}\\
\vdots\notag\\
\varepsilon_i^{(r_i)}=L_{f_i}^{r_i}h_i+L_{g_i}L_{f_i}^{r_i-1}h_iu_i-L_p^{r_i}q. \label{epsilonr}
\end{gather}
The control law $u_i$ can be implemented as
\begin{equation}\label{control1}
u_i=\left(L_{g_i}L_{f_i}^{r_i-1}h_i\right)^{-1}\left(-L_{f_i}^{r_I}h_i+L_p^{r_i}q+v_i\right).
\end{equation}
Then equations (\ref{epsilon1})-(\ref{epsilonr}) result in a $r_i$th-order linear system
\begin{equation}\label{xi}
  \dot{\xi}_i=A_i\xi_i+B_iv_i,
\end{equation}
where $\xi_{ij}\triangleq\phi_{i,j}(x_i)=\varepsilon_{i}^{(j-1)}$ and $\xi_i=col\{\xi_{i1},\xi_{i2},\cdots,\xi_{ir_i}\}$. In order to construct
\begin{equation}
\Phi_i(x_i)=col\{\phi_{i,1},\cdots,\phi_{i,r_i},\phi_{i,r_i+1},\cdots,\phi_{i,n_i}\},
\end{equation}
such that the $i$th follower system can be transformed into (\ref{thm22})(\ref{thm23}), we need to find a group of function $\phi_{i,j},r_i< j\leq n_i$ such that $\Phi_{i*}$ is nonsingular and $L_{g_i}\psi_{i,j}(x_i)=0, r_i< j\leq n_i$. For $1\leq j\leq r_i$, one may notice that
\begin{equation}
d\phi_{i,j}(x_i)=d\varepsilon_{i}^{(j-1)}=\frac{\partial}{\partial x_i^T}(L_{f_i}^{j-1}h_i-dL_p^{j-1}q)=dL_{f_i}^{j-1}h_i.
\end{equation}
Hence,
\begin{align}
&rank\{dh_i,\cdots,dL_{f_i}^{r_i-1}h_i,d\psi_{i,r_i+1},\cdots,d\psi_{i,n_i}\}\notag\\
=&rank\{d\phi_{i,1},\cdots,d\phi_{i,r_i},d\psi_{i,r_i+1},\cdots,d\psi_{i,n_i}\}=n_i.
\end{align}
By setting $\phi_{i,j}=\psi_{i,j},r_i< j\leq n_i$, the tracking error system transformed from the $i$th follower by $(\xi_i^T,\theta^T)^T=\Phi_i(x_i)$ can be expressed as (\ref{thm22})(\ref{thm23}). It indicates that the tracking error system and normal form (\ref{normal}) have the same internal dynamics. Moreover, combining with (\ref{epsilon1})-(\ref{epsilonr}) and (\ref{normal}), we know $\zeta_i=\zeta_0+\xi_i$. Since the stability of zero dynamic $\dot{\theta}_i=\gamma_i(0,\theta_i)$ implies the stability of corresponding internal dynamic (\ref{thm23})\cite{Slotine1991Applied}, the stability of tracking error system can be guaranteed by employing $v_i=K_i\xi_i,$ such that (\ref{thm22}) is stable.
\end{IEEEproof}

Theorem \ref{thm2} gives a purely decentralized control law for the $i$th subsystem. This kind of control law can only be applied to the case where the leader can communicate with all the followers. In the paper, we are supposed to compose the purely decentralized control law (\ref{control1})(\ref{control1v}) and the distributed observer (\ref{dov1})(\ref{dov2}) to further obtain the following distributed control law:
\begin{gather}
\dot{\hat{\eta}}_i=A_0\hat{\eta}_i+a(C\hat{\eta}_i)+cF\sum_{i=1}^Na_{ij}\left(\hat{\eta}_j-\hat{\eta}_i\right)+b_i\left(\hat{\eta}_i-\eta_0\right),\label{control21}\\
\begin{split}
&\hat{u}_i=\left(L_{g_i}L_{f_i\left(x_i,\hat{\eta}_i\right)}^{r_i-1}h_i\right)^{-1}\left(-L_{f_i\left(x_i,\hat{\eta}_i\right)}^{r_I}h_i+L_{p\left(\hat{\eta}_i\right)}^{r_i}q\left(\hat{\eta}_i\right)+v_i\right).\label{control22}
\end{split}
\end{gather}
Despite all this, whether the closed-loop system controlled by a distributed control law based on state estimation is stable is indeed the problem that needs to be further demonstrated.
\begin{theorem}\label{thm3}
The leader-following output tracking problem including leader system (\ref{leadersys1})(\ref{leadersys2}) and follower systems (\ref{followsys1})(\ref{followsys2}) can be solved by distributed control law (\ref{control21})(\ref{control22}) if there exists a distributed observer for leader system. In other words, the distributed control law satisfies certainty equivalence principle.
\end{theorem}

\begin{IEEEproof}
We only need to show the tracking error system convergence to zero under (\ref{control21})(\ref{control22}). For simplifying the symbols, we denote $\upsilon_i=L_{g_i}L_{f_i}^{r_i-1}h_iu_i$, $\hat{\upsilon}_i=L_{g_i}L_{f_i}^{r_i-1}h_i\hat{u}_i$, $\tilde{\upsilon}_i=\hat{\upsilon}_i-\upsilon_i$. By substituting (\ref{control22}) into (\ref{epsilonr}), the $r_i$th derivative of tracking error can be rewritten as
\begin{align*}
\varepsilon_i^{(r_i)}=L_{f_i}^{r_i}h_i+L_{g_i}L_{f_i}^{r_i-1}h_i\left(\hat{u}_i-u_i+u_i\right)-L_p^{r_i}q=v_i+\tilde{\upsilon}_i.
\end{align*}
Then the tracking error system is
\begin{equation}\label{xi2}
  \dot{\xi}_i=A_i\xi_i+B_i\left(v_i+\tilde{\upsilon}_i\right)=(A_i+B_iK_i)\xi_i+B_i\tilde{\upsilon}_i.
\end{equation}
Since $A_i+B_iK_i$ is Hurwize, there exists a positive definite solution $\mathcal{Q}_i$ such that
\begin{equation}
\mathcal{Q}_i(A_i+B_iK_i)+(A_i+B_iK_i)^T\mathcal{Q}_i=-2I_{n_i}.
\end{equation}
Differentiating $V(\xi_i)=1/2\xi_i^T\mathcal{Q}_i\xi_i$ and using Cauchy-Schwarz inequality yields
\begin{align}
\dot{V}(\xi_i)&=-\|\xi_i\|^2+\xi_i^T\mathcal{Q}B_iv_i\notag\\
&\leq -\|\xi_i\|^2+\bar{\sigma}(\mathcal{Q}_i)\|\xi_i\|\|B_i\tilde{\upsilon}_i\|\notag\\
&\leq -\|\xi_i\|^2+\bar{\sigma}(\mathcal{Q}_i)\|\xi_i\|\|\tilde{\upsilon}_i\|\to -\|\xi_i\|^2<0.
\end{align}
It is because $\lim_{t\to\infty}w(t)-\hat{\eta}_i(t)=0$, or equivalently $\lim_{t\to\infty}\tilde{\upsilon}_i=0$, under the condition that the distributed observer exists for leader system. We thus obtain $\lim_{t\to\infty}\xi_i=0$ from $\lim_{t\to\infty}V(\xi_i)<0$.
\end{IEEEproof}

\begin{remark}
According to Section \ref{sec2.4}, the stabilization of zero dynamics means the system is a minimum phase system. Theorem \ref{thm3} can thus be applied to deal with the nonlinear affine system whose relative degree is less than $n$. Therefore,  Theorem 1 in \cite{LiuA} can be regard as a specific case of Theorem \ref{thm3}. The latter will degenerates to the former when each subsystem is supposed to be in a special form
\begin{align*}
\dot{x}_{is}&=x_{(s+1)i}, s=1,2,\cdots,n_i,\\
\dot{x}_{ir}&=f_i(x_i,w)+u_i,\\
y_i&=x_{i1}.
\end{align*}
This system is equivalent to the case where the relative order of each follower is $n_i$.
\end{remark}

\begin{remark}
Theorem \ref{thm3} can be further weakened by using quasi-normal form instead of normal form (\ref{thm22})(\ref{thm23}). Though the normal form always exists for every SISO system, it is difficult to calculate the corresponding diffeomorphism because one need to solve some PDEs.
\end{remark}

\subsection{Distributed observer for MIMO system}

Suppose the dimension of the input and output of the leader system and all follower systems are $m$, i.e., $y_0,y_i,u_i\in\mathbb{R}^m$. The MIMO affine nonlinear dynamics of followers can be described as
\begin{gather}
\dot{x}_i=f_i(x_i,w)+\sum_{j=1}^mg_{ij}(x_j)u_j,\label{mimo1}\\
y_i=\left[h_{i1}(x_i),h_{i2}(x_i),\cdots,h_{im}(x_i)\right]^T.\label{mimo2}
\end{gather}

\begin{definition}
For a MIMO affine nonlinear system (\ref{mimo1})(\ref{mimo2}), suppose $\mathcal{U}$ is neighborhood of a point $x^0$. The (vector) relative degree of this system is $r_{i1},\cdots,r_{im}$ if the following two conditions are fulfilled:\\
(1)\ $L_{g_{ij}}L_{f_i}^lh_{ik}=0$ if for $\forall x\in\mathcal{U}$ and for all $1\leq j\leq m$, $1\leq k\leq m$ and $1\leq l\leq r_{ik}$;\\
(2)\ The following matrix is nonsingular at $x^0$,
\begin{equation*}
\mathcal{A}_i=\begin{pmatrix}L_{g_{i1}}L_{f_i}^{r_{i1}-1}h_{i1}(x^0)&\cdots&L_{g_{im}}L_{f_i}^{r_{i1}-1}h_{i1}(x^0)\\\vdots&\ddots&
\vdots\\L_{g_{i1}}L_{f_i}^{r_{im}-1}h_{im}(x^0)&\cdots&L_{g_{im}}L_{f_i}^{r_{im}-1}h_{im}(x^0)\end{pmatrix}.
\end{equation*}
\end{definition}

For saving of analysis, we still assume that the relative order of each follower is equal everywhere in the whole space, and we further assume that the $i$th follower with relative degree ${r_{i1},r_{i2},\cdots,r_{im}}$ satisfies $\sum_{k=1}^mr_{ik}\leq n_i$. Following the calculation procedure of tracking error system, a coordinate transformation can be defined as:
\begin{gather}
\xi_{ij}^k(x_i)=L_{f_i}^{j-1}h_{ik}(x_i)-L_p^{j-1}q_k\triangleq\left(\varepsilon_i^k\right)^{(j-1)},\label{trans}\\
i=1,2,\cdots,N,\ k=1,2,\cdots,m, j=1,2,\cdots,\ r_{ik}.\notag
\end{gather}
Then the $i$th subsystem can be transformed into $m$ groups equations ($k=1,2,\cdots,m$):
\begin{gather}
\left(\varepsilon_i^k\right)^{(1)}=L_{f_i}h_{ik}(x_i)-L_pq_k, \\
\vdots\ \ \ \notag\\
\left(\varepsilon_i^k\right)^{(r_{ik}-1)}=L_{f_i}^{r_{ik}-1}h_{ik}(x_i)-L_p^{r_{ik}-1}q_k,\\
\left(\varepsilon_i^k\right)^{(r_{ik})}=L_{f_i}^{r_{ik}}h_{ik}(x_i)+\sum_{j=1}^mL_{g_{ij}}L_{f_i}^{r_{ik}-1}h_{ik}u_j-L_p^{r_{ik}}q_k.
\end{gather}
By denoting
\begin{align*}
 &\beta_{ik}(x_i,\omega)=L_{f_i}^{r_{ik}}h_{ik}(x_i)-L_p^{r_{ik}}q_k, \\
 &a_{ikj}(x_i,\omega)=L_{g_{ij}}L_{f_i}^{r_{ik}-1}h_{ik},\\ &\varepsilon_i=col\left\{\left(\varepsilon_i^1\right)^{(r_{i1})},\left(\varepsilon_i^2\right)^{(r_{i2})},\cdots,\left(\varepsilon_i^k\right)^{(r_{im})}\right\},
\end{align*}
the $\left(\varepsilon_i^k\right)^{(r_{ik})}$-dynamic can be introduced in a compact form:
\begin{equation}\label{epsilonm}
\varepsilon_i=\beta_i(x_i,\omega)+\mathcal{A}_iu_i,
\end{equation}
where
\begin{gather*}
  \mathcal{A}_i(x_i,\omega)=\begin{pmatrix}a_{i11}&a_{i12}&\cdots&a_{i1m}\\a_{i21}&a_{i22}&\cdots&a_{i2m}\\\vdots&\vdots&\ddots&\vdots\\a_{im1}&a_{im2}&\cdots&a_{imm}\end{pmatrix},\\
  \beta_i(x_i,\omega)=\begin{pmatrix}\beta_{i1}(x_i,\omega)\\\beta_{i2}(x_i,\omega)\\\vdots\\\beta_{im}(x_i,\omega)\end{pmatrix}, u_i=\begin{pmatrix}u_{i1}\\u_{i2}\\\vdots\\u_{im}\end{pmatrix}.
\end{gather*}

Referring to the definition of relative degree, we know $\mathcal{A}_i$ is invertible. Then a purely decentralized control law for $i$th follower can thus be implemented as
\begin{equation}\label{control3}
u_i=\mathcal{A}_i^{-1}\left(-\beta_i(x_i,\omega)+v_i\right).
\end{equation}
Sequentially, a linear error dynamic can be obtained by combining equations (\ref{trans})(\ref{epsilonm}) and (\ref{control3}):
\begin{equation}
\dot{\xi}_i^k=A_{ik}\xi_i^k+B_{ik}v_{ik},
\end{equation}
where $v_i=col\{v_{i1},v_{i2},\cdots,v_{im}\}$, $\xi_i^k=col\left\{\xi_{ij}^k\right\}_{j=1}^{r_{ik}}$, $\xi_{ij}^k=\left(\varepsilon_i^k\right)^{(j-1)}$, and
\begin{gather*}
A_{ik}=\begin{bmatrix}0&I_{r_{ik}-1}\\0&0\end{bmatrix}\in\mathbb{R}^{r_{ik}\times r_{ik}},\\ b_{ik}=\begin{bmatrix}0&0&\cdots&1\end{bmatrix}^T\in\mathbb{R}^{r_{ik}}.
\end{gather*}
Let $\xi_i=col\left\{\xi_i^k\right\}_{k=1}^m$, $A_i=diag\{A_{i1},A_{i2},\cdots,A_{im}\}$ and $B_i=diag\{B_{i1},B_{i2},\cdots,B_{im}\}$. Then there is a diffeomorphism $\Phi_i$ such that the dynamic of the $i$th follower can be transformed into
\begin{align}
\dot{\xi}_i&=A_i\xi_i+B_iv_i,\label{sys1}\\
\dot{\theta}_i&=\gamma_i(\theta_i,\zeta_i)+\sum_{j=1}^m\rho_{ij}(\theta_i,\zeta_i)u_j.\label{sys2}
\end{align}
According to the knowledge of Theorem \ref{thm2}, we know (\ref{sys2}) is indeed the internal dynamic of \ref{mimo1}, where $\zeta_i=\xi_i+\zeta_0$ with $\zeta_0=col\left\{\zeta_0^k\right\}_{k=1}^m$ and $\zeta_0^k=\left\{L_p^jq(w)\right\}_{j=0}^{r_{ik}-1}$, and the smooth nonlinear function $\gamma_i(\cdot)$ and $\rho_{ij}(\cdot)$ can be obtained following the computation process of quasi-normal form for MIMO affine nonlinear system \cite{Isidori2003Nonlinear}. Note that $\rho_{ij}$ in (\ref{sys2}) could be designed to zero \cite{Isidori2003Nonlinear} if the distribution $\mathcal{D}=span\{g_{i1},g_{i2},\cdots,g_{im}\}$ is involutive. Then the stability of tracking error system (\ref{sys1})(\ref{sys2}) can be ensured by a linear feedback control
\begin{equation}\label{control3v}
v_i=K_i\xi_i,
\end{equation}
if the zero dynamic corresponding to internal dynamic (\ref{sys2}) is stability, where $K_i=diag\{K_{i1},K_{i2},\cdots,K_{im}\}$ is designed to make $A_{ik}+B_{ik}K_{ik}$ be Hurwitz for all $k=1,2,\cdots,m$.

Similar to the previous section, we need to develop the distributed control law corresponding to MIMO system by composing purely decentralized control law (\ref{control3})(\ref{control3v}) and distributed observer (\ref{dov1})(\ref{dov2}) and prove the certainty equivalently principle.

\begin{theorem}\label{thm5}
The leader-following output tracking problem composed of leader system (\ref{leadersys1})(\ref{leadersys2}) and incompletely controllable follower systems (\ref{sys1})(\ref{sys2}) can be solved by distributed control law
\begin{gather}
\dot{\hat{\eta}}_i=A_0\hat{\eta}_i+a(C\hat{\eta}_i)+cF\sum_{i=1}^Na_{ij}\left(\hat{\eta}_j-\hat{\eta}_i\right)+b_i\left(\hat{\eta}_i-\eta_0\right),\label{control41}\\
\hat{u}_i=\hat{\mathcal{A}}_i^{-1}\left(-\hat{\beta}_i(x_i)+v_i\right). \label{control42}
\end{gather}
if there exists a distributed observer for leader system. In (\ref{control42}), $\hat{\mathcal{A}}_i=\mathcal{A}_i(x_i,\hat{\eta}_i)$ and $\hat{\beta}_i=\beta_i(x_i,\hat{\eta}_i)$.
\end{theorem}

\begin{IEEEproof}
Substituting (\ref{control42}) into (\ref{epsilonm}), we have
\begin{equation}
\varepsilon_i=\beta_i(x_i)+\mathcal{A}_i\left(\hat{u}_i-u_i+u_i\right)=v_i+\mathcal{A}_i\left(\hat{u}_i-u_i\right).
\end{equation}
Let $\tilde{\upsilon}_{ikj}=a_{ikj}\left(\hat{u}_{ij}-u_{ij}\right)$ and we can obtain by combining (\ref{epsilonm})
\begin{equation}
\left(\varepsilon_i^k\right)^{(r_{ik})}=v_{ik}+\sum_{j=1}^m\tilde{\upsilon}_{ikj}.
\end{equation}
Then the tracking error system of the $i$th subsystem with the $k$th output is
\begin{equation}\label{xik}
\dot{\xi}_i^k=A_{ik}\xi_i^k+B_{ik}\left(v_{ik}+\sum_{j=1}^m\tilde{\upsilon}_{ikj}\right).
\end{equation}
Noticing that (\ref{xik}) has the same form as (\ref{xi2}). Thus we can prove the solution of satisfies (\ref{xik}) $\xi_i^k(t)\to 0$ if $\lim_{t\to\infty}\tilde{\upsilon}_{ikj}=0$, and the latter can be indicated by (\ref{control41}) directly.
\end{IEEEproof}

\section{Simulation}
\label{sec6}
Firstly, ESSML system is used to show that our novel distributed nonlinear observer based on geometric conditions can be applied to some nonlinear systems that fails to satisfy \cite{LiuA}'s assumption. Then we simulate the distributed observer-based control frame with \emph{Van der Pol} system as leader and an incompletely controllable minimum phase system as two followers. On the one hand, the second example shows that for a nonlinear leader who can satisfy \cite{LiuA}'s hypothesis and geometric conditions, our method can obtain the same distributed observer performance as \cite{LiuA}'s method. On the other hand, our purely decentralized control law based on zero dynamics can make the minimum phase affine nonlinear system which is not completely controllable track the leader's output.
\subsection{Simulation with ESSLM system}
\begin{figure}[!t]
  \centering
  \includegraphics[width=6cm]{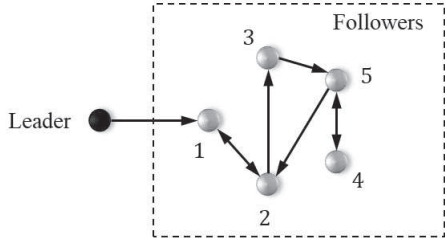}\\
  \caption{A directed communication graph between followers}\label{graph-d}
\end{figure}
\begin{figure}[!t]
  \centering
  \includegraphics[width=8cm]{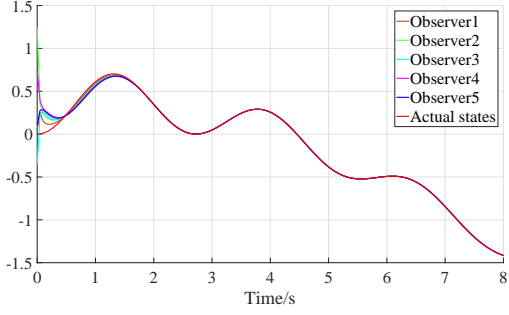}\\
  \caption{State estimate and actual state of $\omega_1$}\label{ESSLM1}
\end{figure}
\begin{figure}[!t]
  \centering
  \includegraphics[width=8cm]{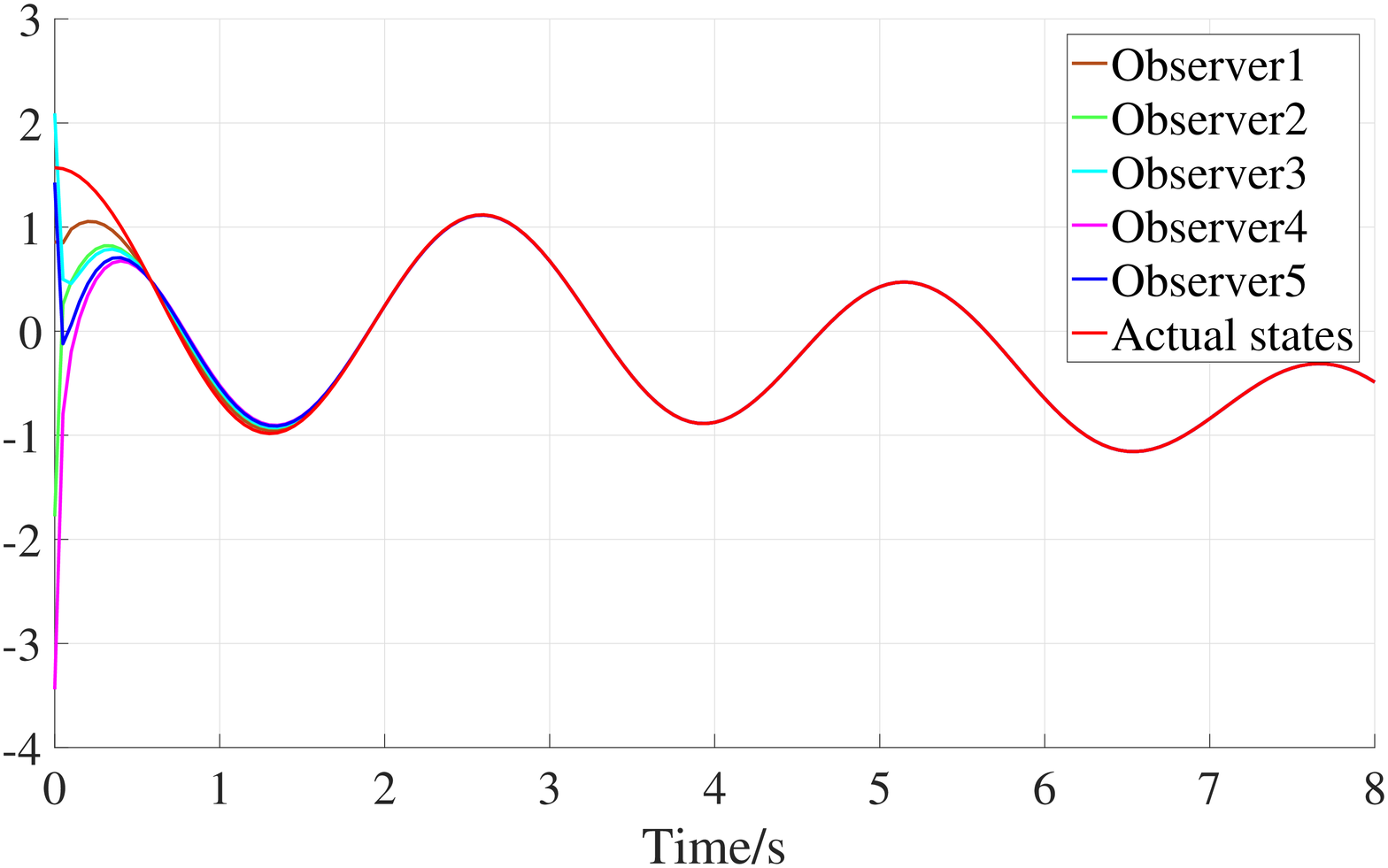}\\
  \caption{State estimate and actual state of $\omega_2$}\label{ESSLM2}
\end{figure}
\begin{figure}[!t]
  \centering
  \includegraphics[width=8cm]{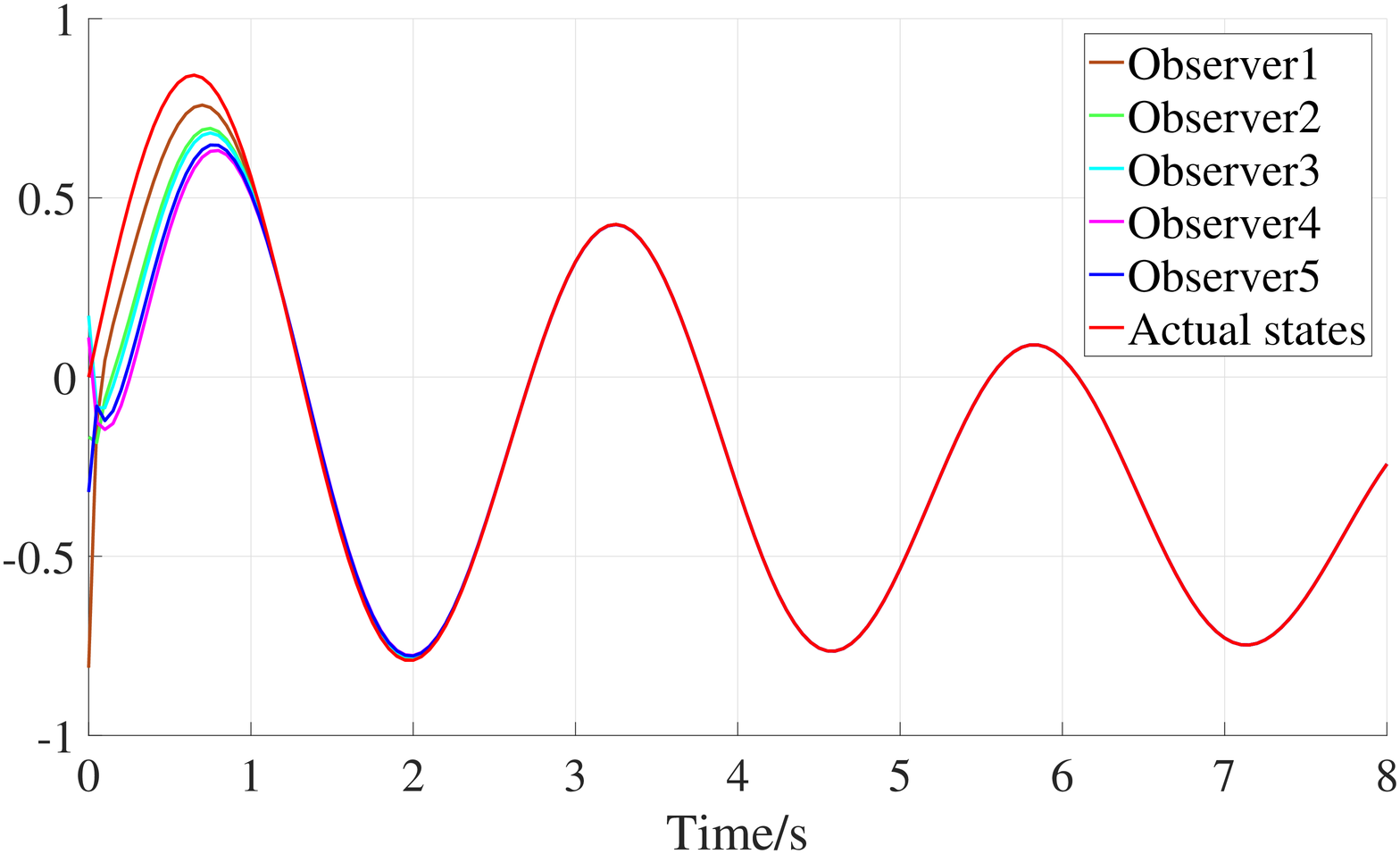}\\
  \caption{State estimate and actual state of $\omega_3$}\label{ESSLM3}
\end{figure}
\begin{figure}[!t]
  \centering
  \includegraphics[width=8cm]{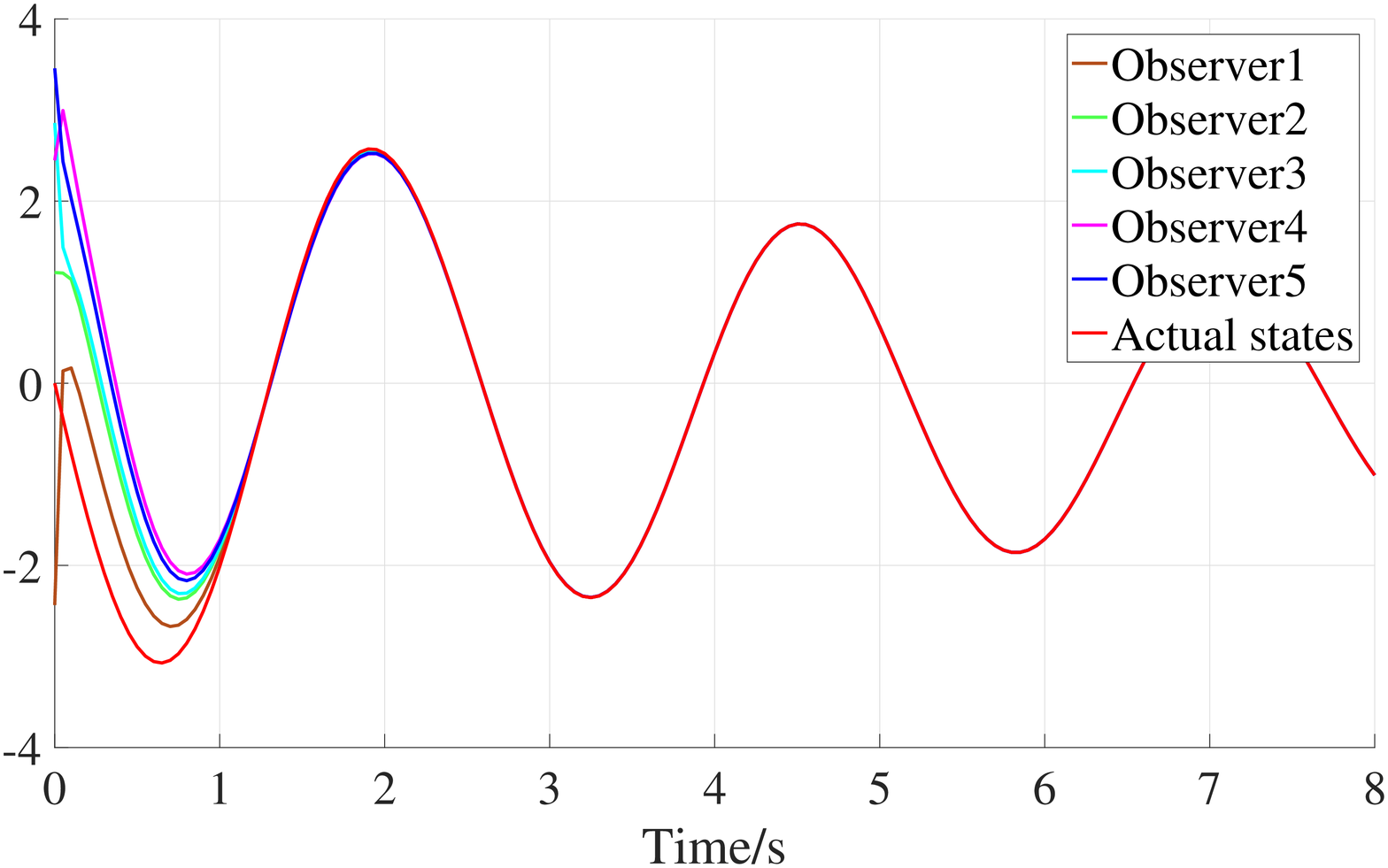}\\
  \caption{State estimate and actual state of $\omega_4$}\label{ESSLM4}
\end{figure}

Section \ref{sec3.2} has proved that ESSLM system (\ref{esslm}) satisfies geometric conditions in Lemma \ref{lemma1}. Suppose there are five followers and the communication graph between leader and followers is showed in Figure \ref{graph-d}. In (\ref{esslm}), we set the length of the Linkage be $ 2d=0.2{\rm m}$, the mass of Linkage be $m=1{\rm kg}$, the Rotational
inertia be $J_1=5{\rm kg\cdot m^2}$, $J_2=2{\rm kg\cdot m^2}$, the Viscosity friction coefficient be $F_1=0.5,F_2=0.55$, and the torsional elastic coefficient of elastic shaft be $\mathcal{K}={\rm Nm/rad}$. The acceleration of gravity is approximately taken as $g=10{\rm m/s^2}$. From the calculation in Section \ref{sec3.2}, we can obtain a diffeomorphism $\eta_0=\Phi(\omega)$ such that
\begin{equation}
\eta_0=\begin{bmatrix}
  0.33&      0.244& 3.33&  0.889\\
 3.33&       0.916&    0& 0.1\\
    0&       0.375&    0&    1\\
    0&         1&    0&    0
\end{bmatrix}\omega.
\end{equation}
and
\begin{equation}
\omega=\begin{bmatrix}
   0&   0.3&    -0.03&  -0.264\\
    0&      0&         0&          1\\
0.3& -0.03& -0.264&  0.053\\
    0&      0&         1&       -0.375
\end{bmatrix}\eta_0.
\end{equation}
By calculating $\left.\frac{\partial\Phi}{\partial\omega^T}p(\omega)\right|_{\omega=\Phi^{-1}(\eta_0)}$, we have the observer canonical form of ESSLM:
\begin{align*}
\dot{\eta}_0&=\begin{bmatrix}\dot{\eta}_{01}\\\dot{\eta}_{02}\\\dot{\eta}_{03}\\\dot{\eta}_{04}\end{bmatrix}
=\begin{bmatrix}0&0&0&0\\1&0&0&0\\0&1&0&0\\0&0&1&0\end{bmatrix}\eta_{0}+\begin{bmatrix}- \frac{4}{9}\cos{\eta_{04}}\\
- \frac{67}{90}\eta_{04} - \frac{1}{20}\cos{\eta_{04}}\\
- \frac{21299}{3600}\eta_{04} - \frac{1}{2}\cos{\eta_{04}}\\
- \frac{3}{8}\eta_{04}
\end{bmatrix},\\
y_0&=\eta_{04}.
\end{align*}
Therefore, a distributed observer for this leader can be designed by (\ref{dov1})(\ref{dov2}) with $c=5$ and
\begin{equation*}
F=\begin{bmatrix}
    5.07  &  0.50 &   0.02 &   0.00\\
    0.50  &  5.05 &   0.50 &   0.02\\
    0.02  &  0.50 &   5.05 &   0.50\\
    0.00  &  0.02 &   0.50 &   4.98
\end{bmatrix}.
\end{equation*}
The initial states of original system is performed with $[0,\pi/2,0,0]^T$ and that of each agent are generated randomly. Figures \ref{ESSLM1} - Figure \ref{ESSLM4} show the comparison between the actual states of the leader system and the state estimates generated by each local observer. These figures illustrate that the state estimates of the distributed observer converge quickly to the actual states, which verifies the effectiveness of our new method. In other words, we can design a distributed observer for ESSLM system, a leader system which cannot be handled by \cite{LiuA}, and obtain excellent dynamic performance.

\begin{figure}[!t]
  \centering
  \includegraphics[width=8cm]{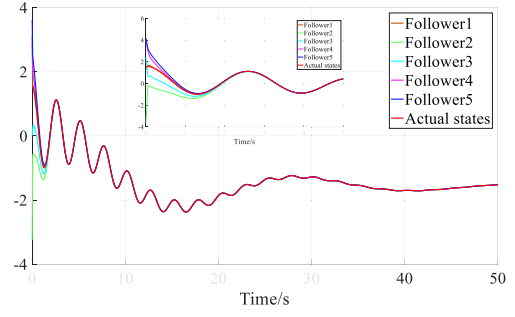}\\
  \caption{Phase of linkage}\label{ESSLMout}
\end{figure}

\begin{figure}[!t]
  \centering
  \includegraphics[width=8cm]{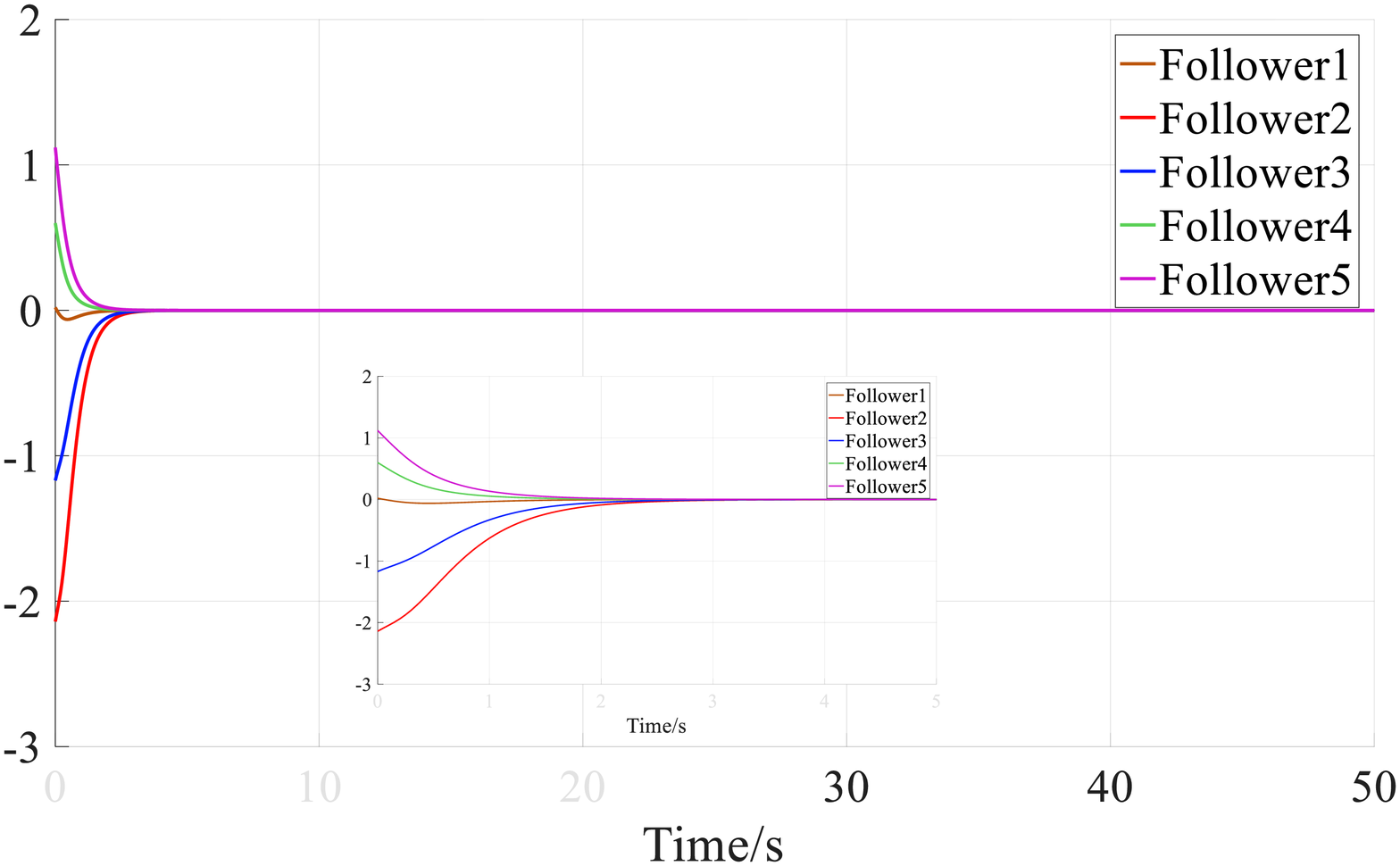}\\
  \caption{Tracking error of followers}\label{ESSLMerror}
\end{figure}

Five followers ($i=1,\cdots,5$) are also chosen as ESSLM system:
\begin{equation}\label{esslm-follow}
\begin{split}
\begin{bmatrix}\dot{\omega}_{i1}\\ \dot{\omega}_{i2}\\ \dot{\omega}_{i3}\\ \dot{\omega}_{i4}\end{bmatrix}
&=\begin{bmatrix}\omega_{i3}\\ \omega_{i4}\\
\frac{4}{3}\omega_{i1}-\frac{8}{9}\omega_{i2}-\frac{1}{10}\omega_{i3}\\
\frac{10}{3}\omega_{i1}-5\omega_{i2}-\frac{11}{40}\cos{\omega_{i2}}-\frac{F_2}{J_2}\omega_{i4}\end{bmatrix}
+\begin{bmatrix}0\\0\\ \frac{1}{5}\\0\end{bmatrix}u_i,\\
y_i&=q(\omega)=\omega_{i2}.
\end{split}
\end{equation}
Unlike the leader system (\ref{esslm}), follower system contains an affine nonlinear control input. Note that ESSLM system is a completely controllable system, thus it can be controlled be feedback linearization directly. The initial value of the follower systems are chosen randomly. Figure \ref{ESSLMout} illustrates that the output of all followers can track the output of leader. Figure \ref{ESSLMerror} demonstrates the tracking error dynamic between leader and followers.

\subsection{Simulation with Van der Pol system}
Suppose the leader obey the follow \emph{Van der Pol} system
\begin{align}\label{vdp}
\dot{w}_1&=w_2,\notag\\
\dot{w}_2&=-w_1+(1-w_1^2)w_2,\\
y_0&=w_1.\notag
\end{align}
The communication graph between all followers is designed as Figure \ref{graph-d}. It is easy to check $[dq(w),dL_pq(w)]=I_2$, i.e., \emph{Van der Pol} system satisfies observability condition (Lemma \ref{lemma1} (1)). Utilizing Lemma \ref{lemma1}, we can also calculate $\tau(w)=[0,1]^T$, hence, $[\tau(w),ad_p\tau(w)]=0$. Thus (\ref{le1})(\ref{le2}) are satisfied. So we can find a coordinate transformation (The method of solve this coordinate transformation refer to appendix \cite{Xia1989Nonlinear})
\begin{equation*}
  \begin{bmatrix}\eta_{01}\\\eta_{02}\end{bmatrix}
  =\Phi(w)
  =\begin{bmatrix}-w_1+\frac{1}{3}w_1^3+w_2\\w_1\end{bmatrix},
\end{equation*}
and its inverse information
\begin{equation*}
  \begin{bmatrix}w_1\\w_2\end{bmatrix}
  =\Phi^{-1}(\eta_0)
  =\begin{bmatrix}\eta_{02}\\\eta_{01}+\eta_{02}-\frac{1}{3}\eta_{02}^3\end{bmatrix},
\end{equation*}
such that leader system (\ref{vdp}) is transformed in observable canonical form
\begin{align*}
 \dot{\eta}_0&=A_0\eta_0+a(y_0),\\
 y_0&=\eta_{02},
\end{align*}
where,
\begin{equation*}
  A_0=\begin{bmatrix}0&0\\1&0\end{bmatrix},\ \ \ a(y_0)=\begin{bmatrix}-\eta_{02}\\\eta_{02}-\frac{1}{3}\eta_{02}^3\end{bmatrix}.
\end{equation*}

\begin{figure}[!t]
  \centering
  \includegraphics[width=8cm]{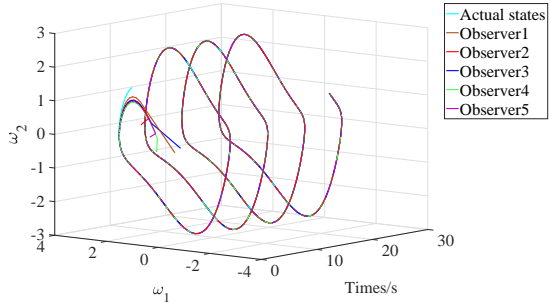}\\
  \caption{Profile on the phase portraits of the leader and the distributed observer with geometric conditions.}\label{vdp1}
\end{figure}

\begin{figure}[!t]
  \centering
  \includegraphics[width=8cm]{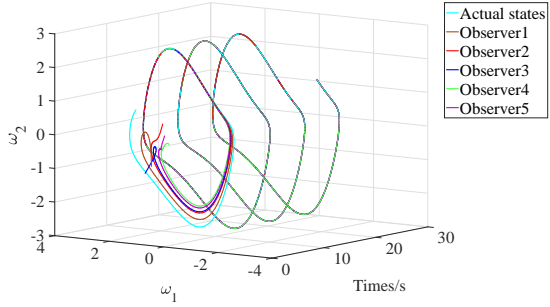}\\
  \caption{Profile on the phase portraits of the leader and the distributed observer with Taylor conditions.}\label{vdp2}
\end{figure}

Let $c=10$ so that it satisfies conditions (\ref{c1})(\ref{c2}). Figure \ref{vdp1} demonstrates that all the state estimate generated by followers can converge to the actual states. These simulation results indicate the correction of Theorem \ref{thm1}. Furthermore, comparing to Figure \ref{vdp2}, the distributed observer obtained by \cite{LiuA}'s model, it can be seen that our novel distributed observer has faster convergence speed under the same coupling gain.

\begin{figure}[!t]
  \centering
  \includegraphics[width=8cm]{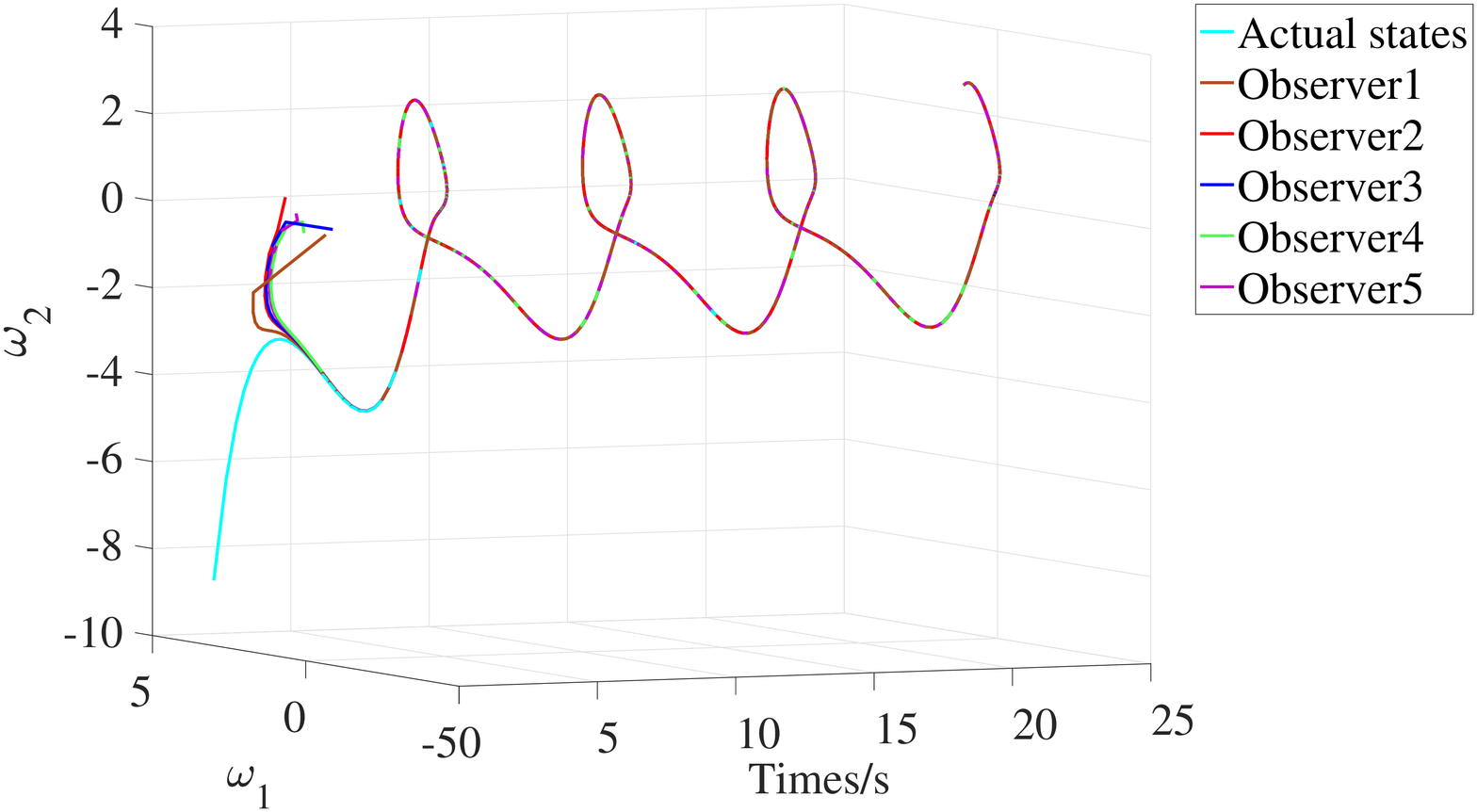}\\
  \caption{The method based on geometric conditions proves that the Van der Pol system can actually have a globally convergent distributed observer}\label{vdp3}
\end{figure}

In addition, \cite{LiuA}'s assumption limits the application scope of their distributed observer to a compact set containing the origin. For example, for van der Pol system, the initial value of their leader needs to be selected in $\|\eta_0(0)\|\leq2\sqrt{2}$. Actually, the distributed observer designed for Van der Pol system can be globally convergent owing to it meets geometric conditions globally. Figure \ref{vdp3} shows the convergence performance when $\|\eta_0(0)\|$ is chosen outside of $\|\eta_0(0)\|\leq2\sqrt{2}$.
\begin{figure}[!t]
  \centering
  \includegraphics[width=8cm]{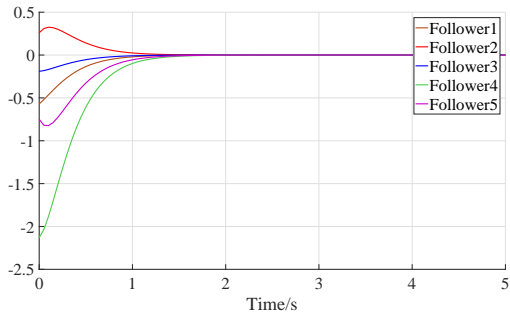}\\
  \caption{Tracking error of followers}\label{output}
\end{figure}

Assume followers 1,3,5 satisfy a nonlinear system
\begin{equation}
\begin{split}
\dot{x}_{i1}&=x_{i1}+x_{i2},\\
\dot{x}_{i2}&=x_{i1}x_{i2}^{a_i}+u_i,\\
y_i&=x_{i1},
\end{split}
\end{equation}
where $a_i$ for all $i=1,2,\cdots,N$ are parameters depended on $i$th subsystem. One can check that every subsystem has relative degree $2$ under the given output $y_i$. Thus there is a coordinate transformation
\begin{align*}
\xi_{i1}&=x_{i1}-w_1,\\
\xi_{i2}&=x_{i1}+x_{i2}-L_pq(w),\\
y_i&=\xi_{i1}+w_1.
\end{align*}
Then the purely decentralized control law can be designed as
\begin{equation}\label{controllaw1}
u_i=-x_{i1}-x_{i2}+L_pq(\hat{\eta}_i)+v_i.
\end{equation}
On the other hand, followers 2 and 4 are in the form of \cite{Slotine1991Applied}
\begin{equation}\label{f2}
\begin{split}
\dot{x}_{i1}&=-x_{i1}+e^2x_{i2}u_i,\\
\dot{x}_{i2}&=2x_{i1}x_{i2}+\sin x_{i2}+\frac{1}{2}u_i,\\
\dot{x}_{i3}&=2x_{i2},\\
y_i&=x_{i3}.
\end{split}
\end{equation}
This system is not incompletely controllable, hence, it cannot be controlled by \cite{LiuA}'s purely decentralized control law. By transforming it into quasi-normal form:
\begin{align}
\dot{\xi}_{i1}&=\xi_{i2},\notag\\
\dot{\xi}_{i2}&=2\left(-1+\theta_i+e^{\xi_{i2}}\right)\xi_{i2}+2\sin \frac{\xi_{i2}}{2}-L_pq(w)+u_i,\notag\\
\dot{\theta}_i&=\left(1-\theta_i-e^{\xi_{i2}}\right)\left(1+2\xi_{i2}e^{\xi_{i2}}\right)-2\sin \frac{\xi_{i2}}{2}e^{\xi_{i2}},\label{interdy}\\
y_i&=\xi_{i1}+w_1.
\end{align}
we obtain its inter dynamic (\ref{interdy}). Then we can further get zero dynamic by setting $\xi_{i1}=\xi_{i2}=0$:
\begin{equation}
\dot{\theta}_i=-\theta_i.
\end{equation}
It is obviously that the zero dynamic of (\ref{f2}) is stable. Hence, we can design the purely decentralized control law for this system:
\begin{equation}\label{controllaw2}
u_i=-2\left(-1+\theta_i+e^{\xi_{i2}}\right)\xi_{i2}-2\sin \frac{\xi_{i2}}{2}+L_pq(w)+v_i.
\end{equation}

Then the distributed control law of this leader-following problem can be constructed by replacing state estimates $\Phi^{-1}\hat{\eta}_i$ generated by distributed observer with $\omega$ in purely decentralized control law (\ref{controllaw1}) and (\ref{controllaw2}). The initial states of each subsystem are chosen randomly and the pole of the feedback linearization system is allocated at $-2,-6$. Figure \ref{output} shows the tracking error of subsystems to external system under the distributed control law. It can be seen that the leader-following consensus is achieved.

\section{Conclusion}
\label{sec7}
This paper has proposed a novel distributed nonlinear observer based on geometric conditions. Within this method, a special assumption on leader system constrained by \cite{LiuA} has been replaced with a group of geometric conditions. As a result, our distributed nonlinear observer can be applied for some nonlinear system which fails to fulfill \cite{LiuA}'s assumption, such as ESSLM system and most of first-order nonlinear system. We have proved that our distributed nonlinear observer has an exponentially stable error dynamics for all the output bounded nonlinear system met geometric conditions. Two lemmas corresponding to the spectrum of the matrices are proved as a pioneer to complete the proof. Furthermore, we have developed purely decentralized control law based on zero dynamic proposed in differential geometry. With this advancement, the followers can be chosen as an arbitrary minimum phase affine nonlinear system. The certainty equivalence principle for the distributed observer-based control law including novel distributed nonlinear observer and improved purely decentralized control law has also been prove. ESSLM system and Van der Pol system have been used to simulate our method.


%

\section*{Appendix}
\subsection*{Computation procedure of observable canonical form}
\noindent\textbf{STEP 1} Calculate $\tau_1,\tau_2,\cdots,\tau_r$ of (\ref{le1})(\ref{le2}).

\noindent\textbf{STEP 2} Compute a matrix as follow:
\begin{align}
\tilde{\mathcal{Q}}_i&=\left[\tau_i,ad_{-p}\tau_i,\cdots,ad_{-p}\tau_i\right].\ t=1,2,\cdots,r.\\
\tilde{\mathcal{Q}}&=\left[\tilde{\mathcal{Q}}_1,\tilde{\mathcal{Q}}_2,\cdots,\tilde{\mathcal{Q}}_r\right].
\end{align}

\noindent\textbf{STEP 3} Compute $b_i(y_0)$ defined as
\begin{equation}
b_i(y_0)=\tilde{\mathcal{Q}}^{-1}(x)ad_{-p}^{k_i}\tau_i.
\end{equation}

\noindent\textbf{STEP 4} Solving the following equations:
\begin{equation}
\frac{\partial a(y_0)}{\partial y_{0i}}=b_i(y_0).
\end{equation}
\noindent\textbf{STEP 5} Denote $\pi_i=\sum_{j=1}^ik_j$, then the coordinate transformation $\eta_0=\Phi(x)$ can be calculated by
\begin{equation}
\eta_{0i}(x)=\begin{cases}h_i(x), &\text{if}\ i\in\{\pi_1,\pi_2,\cdots,\pi_r\};\\
L_f\eta_{0,i+1}(x)-a_i(h(x)), &\text{if}\ i\notin\{\pi_1,\pi_2,\cdots,\pi_r\};\end{cases}
\end{equation}
where, $i=1,2,\cdots,s$.

\ifCLASSOPTIONcaptionsoff
  \newpage
\fi



%




\bibliographystyle{IEEEtran}
\bibliography{tex}

\end{document}